\documentclass[a4paper,12pt]{article}
\usepackage{amsfonts}
\usepackage{amssymb,amsmath,amsthm}
\usepackage[matrix,arrow]{xy}

\newtheorem{prop}{Proposition}

\newcommand{\id}{\mathop{\mathrm{id}}\nolimits}

\newcommand{\ctg}{\mathop{\mathrm{ctg}}\nolimits}

\newcommand{\cth}{\mathop{\mathrm{cth}}\nolimits}

\newcommand{\End}{\mathop{\mathrm{End}}\nolimits}

\newcommand{\Aut}{\mathop{\mathrm{Aut}}\nolimits}

\newcommand{\tr}{\mathop{\mathrm{tr}}\nolimits}

\renewcommand{\Im}{\mathop{\mathrm{Im}}\nolimits}
\renewcommand{\Re}{\mathop{\mathrm{Re}}\nolimits}

\numberwithin{equation}{section}
\newcommand*{\hm}[1]{#1\nobreak\discretionary{}%
            {\hbox{$\mathsurround=0pt #1$}}{}}

\newcommand{\A}{{\cal A}}

\renewcommand{\L}{L}
\newcommand{\LL}{{\cal L}}
\newcommand{\rr}{{\mathfrak r}}
\newcommand{\G}{G_{\lambda}^}
\newcommand{\Gg}{{\cal G}_{\lambda}^}

\newcommand{\ph}{\varphi^}
\newcommand{\Ph}{\Phi^}
\newcommand{\e}{e_{\lambda}}
\newcommand{\f}{f_{\lambda}}
\newcommand{\la}{\big\langle}
\newcommand{\ra}{\big\rangle}
\newcommand{\La}{\Big\langle}
\newcommand{\Ra}{\Big\rangle}

\newcommand{\bs}[1]{\backslash\{#1\}}
\newcommand{\lfK}{{\cal K}}

\newcommand{\bU}{\textbf{U}}
\newcommand{\Cyl}{\mathrm{Cyl}}
\newcommand{\slt}{\mathfrak{sl}_2}
\newcommand{\cEF}{\mathfrak{e}_\tau(\widehat{\mathfrak{sl}}_2)}
\newcommand{\utau}{\mathfrak{u}_\tau(\widehat{\mathfrak{sl}}_2)}
\newcommand{\vpint}{\mathop{\hspace{-1pt}{-}\hspace{-12pt}\int}}
\newcommand{\vpints}{\mathop{\hspace{-1pt}{-}\hspace{-10pt}\int}}

\def\fK{\mathfrak{K}}
\def\CC{\mathbb{C}}
\def\frg{\mathfrak{g}}
\def\hg{\hat{\mathfrak{g}}}
\def\r#1{(\ref{#1})}
\def\hx{\hat{x}}

\def\nn{\nonumber}

\newcounter{bbcount}[subsection]
\newcommand{\bb}[1]{\addtocounter{bbcount}{1}{\bf {\thesection}.\arabic{subsection}.\arabic{bbcount}.}
{\bfseries #1}}

\begin{document}

\begin{center}

\hfill ITEP-TH-31/07\\
\hfill ESI-1950\\
\hfill math.QA/yymmnnn\\
\bigskip
{\Large\bf Classical elliptic current algebras}
\par\bigskip\medskip
{\bf Stanislav Pakuliak$^{\bullet\star}$\footnote{E-mail:
pakuliak@theor.jinr.ru}, Vladimir
Rubtsov$^{\circ\star}$\footnote{E-mail:
Volodya.Roubtsov@univ-angers.fr} and Alexey
Silantyev$^{\bullet\circ}$\footnote{E-mail: silant@theor.jinr.ru, silant@tonton.univ-angers.fr}
}\ \
\par\bigskip\medskip
$^\star${\it Institute of Theoretical \& Experimental Physics, 117259
Moscow, Russia}
\par\smallskip
$^\bullet${\it Laboratory of Theoretical Physics, JINR,
141980 Dubna, Moscow reg., Russia}\par\smallskip
$^\circ${\it D\'epartment de Math\'ematiques,
Universit\'e d'Angers,\\
2 Bd. Lavoisier, 49045 Angers, France}
\\
\bigskip
\bigskip
\emph{In memory of Leonid Vaksman}
\\
\bigskip
\bigskip

\end{center}

\thispagestyle{empty}

\begin{abstract}
In this paper we discuss classical elliptic current algebras and show that there are two different choices of commutative {\em test function algebras} on a complex torus leading to two different elliptic current algebras. Quantization of these classical current algebras give rise to two classes of quantized dynamical quasi-Hopf current algebras studied by Enriquez-Felder-Rubtsov and Arnaudon-Buffenoir-Ragoucy-Roche-Jimbo-Konno-Odake-Shi\-raishi. Different degenerations of the classical elliptic algebras are considered. They yield different versions of rational and trigonometric current algebras. We also review the averaging method of Faddeev-Reshetikhin, which allows to restore elliptic algebras from the trigonometric ones.
\end{abstract}

\section{Introduction}
\label{sec1}

Classical elliptic algebras are `quasi-classical limits' of
quantum algebras whose structure is defined by an elliptic
$R$-matrix. The first elliptic $R$-matrix appeared as a matrix
of Boltzmann weights for the eight-vertex model~\cite{Bax1}. This
matrix satisfies the Yang-Baxter equation using which one proves
integrability of the model. An investigation of the eight-vertex
model~\cite{Bax2} uncovered its relation  to the so-called
generalized ice-type model -- the Solid-On-Solid (SOS) model. This is
a face type model with Boltzmann weights which form a matrix
satisfying a {\em dynamical} Yang-Baxter equation.

In this paper we restrict our attention to {\em classical current
algebras} (algebras which can be described by a collection of currents)
related to the classical $r$-matrix and which are quasi-classical limits of
SOS-type {\em quantized elliptic current algebras}. The latter
were introduced by Felder~\cite{F} and the corresponding $R$-matrix
is called usually a {\em Felder $R$-matrix}. In \emph{loc. cit.}
the current algebras were defined by dynamical $RLL$-relations. At
the same time Enriquez and one of authors (V.R.) developed  a
theory of quantum current algebras related to arbitrary genus
complex curves (in particular to an elliptic curve) as a
quantization of certain (twisted) Manin pairs~\cite{ER97}
using {\em Drinfeld's new realization} of quantized current
algebras. Further, it was shown in~\cite{EF98}
that the Felder algebra can be obtained by twisting of the
Enriquez-Rubtsov elliptic algebra. This twisted algebra will be
denoted by $E_{\tau,\eta}$ and it is a quasi-Hopf algebra.

Originally, the dynamical Yang-Baxter equation appeared in ~\cite{GN,F94}. The
fact that elliptic algebras could be obtained as quasi-Hopf
deformations of Hopf algebras was noted first in a special
case in ~\cite{BBB} and was discussed in ~\cite{Fr97}. The full potential
of this idea was realized in papers ~\cite{ABRR97} and
~\cite{JKOS97}. It was explained in these papers how to obtain
the universal dynamical Yang-Baxter equation for the twisted elliptic universal
$R$-matrix from the Yang-Baxter equation for the universal $R$-matrix of the quantum
affine algebra $U_q(\hat{\mathfrak g})$. It was also shown that the image of the twisted
$R$-matrix in finite-dimensional representations coincides with
SOS type $R$-matrix.

 Konno proposed in~\cite{K98} an RSOS type elliptic current algebra (which will be denoted by
$U_{p,q}(\hat{\mathfrak{sl}}_2)$) generalizing some ideas of~\cite{KLP98}. This algebra was studied in detail in ~\cite{JKOS} where it was shown that commutation relations for $U_{p,q}(\hat{\mathfrak{sl}}_2)$ expressed in terms of $L$-operators coincide with the commutation relations of the Enriquez-Felder-Rubtsov algebra up to a shift of the elliptic module by the central element. In the case of the algebra $E_{\tau,\eta}$ the elliptic module is fixed, while in the case of $U_{p,q}(\hat{\mathfrak{sl}}_2),\ p=e^{i\pi\tau},$ it turn out to be a dynamical parameter shifted by the central element. Commutation relations for these algebras coincide when the central charge is zero, but the algebras themselves are different. Furthermore, the difference between these two algebras was interpreted in~\cite{EPR} as a difference in definitions of half-currents (or Gauss coordinates) in $L$-operator representation. The roots of this difference are related to different decomposition types of so-called {\em Green kernels} introduced in~\cite{ER97} for quantization of Manin pairs: they are expanded into Taylor series in the case of the algebra $E_{\tau,\eta}$ and into Fourier series for
$U_{p,q}(\hat{\mathfrak{sl}}_2)$~\cite{JKOS}.

Here, we continue a comparative study of different elliptic current
algebras. Since the Green kernel is the same in both the classical
and the quantum case we restrict ourselves only to the classical
case for the sake of  simplicity. The classical limits of
quasi-Hopf algebras $E_{\tau,\eta}$ and
$U_{p,q}(\hat{\mathfrak{sl}}_2)$ are {\em quasi-Lie bialgebras}
denoted by $\cEF$ and $\utau$ respectively. We will give an
`analytic' description of these algebras in terms of
distributions. Then, the different expansions of Green kernels
will be interpreted as the distributions acting on different test
function algebras. We will call them {\em Green distributions}.
The scalar products for test function algebras which define their
embedding in the corresponding space of distributions are defined by
integration over different contours on the surface.

Let us describe briefly the structure of the paper. Section 2
contains some basic notions and constructions which are used
throughout the paper. Here, we remind some
definitions from ~\cite{ER97}. Namely, we define test function
algebras on a complex curve $\Sigma$, a continuous non-degenerate
scalar product, distributions on the test functions and a
generalized notion of Drinfeld currents associated with these
algebras and with a (possibly infinite-dimensional) Lie algebra
$\mathfrak g$. Hence, our currents will be certain $\mathfrak
g$-valued distributions. Then we review the case when $\mathfrak
g$ is a loop algebra generated by a semi-simple Lie algebra
$\mathfrak a$. We also discuss a centrally and co-centrally
extended version of $\mathfrak g$ and different bialgebra
structures. The latter are based on the notion of Green
distributions and related half-currents.

We describe in details two different classical elliptic current
algebras which correspond to two different choices of the basic
test function algebras (in fact they correspond to two different
coverings of the underlying elliptic curve).

Section 3 is devoted to the construction and comparison of
classical elliptic  algebras $\cEF$ and $\utau$. In the first two
subsections we define elliptic Green distributions for both test
function algebras. We pay special attention to their properties
because they manifest the main differences between the
corresponding elliptic algebras. Further, we describe these
classical elliptic algebras in terms of the half-currents
constructed using the Green distributions. We see how the
half-currents inherit the properties of Green distributions. In
the last subsection we show that the half-currents describe the
corresponding bialgebra structure. Namely, we recall the universal
classical $r$-matrices for both elliptic classical algebras $\cEF$
and $\utau$ and make explicit their relation to the $L$-operators.
Then, the corresponding co-brackets for half-currents are expressed
in a matrix form via the $L$-operators.

 In Section 4 we describe different degenerations of the
classical elliptic current algebras in terms of
degenerations of Green distributions entering the $r$-matrix. The
degenerate Green functions define the $rLL$-relations, the
bialgebra structure and the analytic structure of half-currents.
We do not write out explicitly the bialgebra structure related to
the half-currents of the second classical elliptic algebra: it can
be reconstructed along the lines of subsection~\ref{subsec34}.

We discuss the inverse (in some sense) problem in Section 5. A
way to present the trigonometric and elliptic solutions of
a Classical Yang-Baxter Equation (CYBE) by averaging of the
rational ones was introduced in ~\cite{RF}. Faddeev and Reshetikhin applied the
averaging method to a description of corresponding algebras. Here, we only
represent the elliptic $r$-matrix $\rr^+_\lambda(u-v)$ and
the trigonometric $r$-matrix $\rr^{(c)+}(u-v)$ as an average of
the trigonometric $r$-matrix $\rr^{(b)+}(u-v)$ and the rational $r$-matrix
$\rr^{(a)+}(u-v)$ respectively, for some domains of parameters.

Finally, in the Appendix, we have collected technical and `folklore' definitions, results concerning the
test and distribution algebras on Riemann surfaces. Although some of the results
can be extracted from standard textbooks~\cite{GSh1}, \cite{Vl}, we were not able to
find them in the literature in the form suited for our goals and we have decided to keep
them for the sake of completeness.

\section{Currents and half-currents}

Current realization of the quantum affine algebras and Yangians
was introduced by Drinfeld in~\cite{D88}. In these cases they can
be understood as elements of the space $\A[[z,z^{-1}]]$, where
$\A$ is a corresponding algebra. Here we introduce a more general
notion of currents suitable even for the case when the currents
are expressed by formal integrals instead of formal series (as in
sec.~\ref{sec4}).

\bb{Test function algebras.} Let $\fK$ be a function algebra on a
one-dimensional complex manifold $\Sigma$ with a point-wise
multiplication and a continuous invariant (non-degenerate) scalar
product $\la\cdot,\cdot\ra\colon\fK\times\fK\to\CC$. We shall call
the pair $(\fK,\la\cdot,\cdot\ra)$ {\em a test function algebra}.
The non-degeneracy of the scalar product implies that the algebra
$\fK$ can be extended to a space $\fK'$ of linear continuous
functionals on $\fK$. We use the notation $\la a(u),s(u)\ra$ or
$\la a(u),s(u)\ra_u$ for the action of the distribution
$a(u)\in\fK'$ on a test function $s(u)\in\fK$. Let
$\{\epsilon^i(u)\}$ and $\{\epsilon_i(u)\}$ be  dual bases of
$\fK$. A typical example of the element from $\fK'$ is the series
$\delta(u,z)=\sum_i\epsilon^i(u)\epsilon_i(z)$. This is a
delta-function distribution on $\fK$ because it satisfies $\la
\delta(u,z),s(u)\ra_u=s(z)$ for any test function $s(u)\in\fK$.

\bb{Currents.} Consider an infinite-dimensional complex Lie algebra $\frg$ and an operator
$\hx\colon\fK\to\frg$. The expression
$x(u)=\sum_i \epsilon^i(u)\hx[\epsilon_i]$ 
does not depend on a choice of dual bases in $\fK$ and is called a current corresponding
to the operator $\hx$ ($\hx[\epsilon_i]$ means an action of $\hx$ on $\epsilon_i$).
We should interpret the current $x(u)$ as a $\frg$-valued distribution such that
$\la x(u),s(u)\ra=\hx[s]$. 
That is the current $x(u)$ can be regarded as a kernel of the operator $\hx$ and the
latter formula
 gives its invariant definition.

\bb{Loop algebras} Let $\{\hat x_k\},\ k=1,\ldots, n$ be a finite
number of operators $\hat x_k\colon\fK\to\frg$, where $\frg$ is an
infinite-dimensional space spanned by $\hat x_k[s]$, $s\in\fK$.
Consider the corresponding currents $x_k(u)$. For these currents
we impose the standard commutation relations
\begin{align}\label{cu-com-rel}
[x_k(u),x_l(v)]=C^m_{kl}x_m(u)\delta(u,v)\,
\end{align}
where $C^m_{kl}$ are structure constants of some semi-simple Lie
algebra $\mathfrak a,\ \dim\mathfrak a = n$
(equality~\r{cu-com-rel} is understood in sense of distributions).
These commutation relations equip $\frg$ with a Lie algebra
structure. The Lie algebra $\frg$ defined in such a way can be
viewed as a Lie algebra $\mathfrak a\otimes\fK$ with the brackets
$[x\otimes s(z),y\otimes t(z)]=[x,y]_{\mathfrak a}\otimes
s(z)t(z)$, where $x,y\in\mathfrak a$, $s,t\in\fK$. This algebra
possesses an invariant scalar product
$\la x\otimes s,y\otimes t\ra=(x,y)\la s(u),t(u)\ra_u$, 
where $(\cdot,\cdot)$ an invariant scalar product on $\mathfrak a$ proportional to the Killing form.

\bb{Central extension.} The algebra $\frg=\mathfrak a\otimes\fK$
can be extended by introducing a central element $c$ and a
co-central element $d$. Let us consider the space $\hg=(\mathfrak
a\otimes\fK)\oplus\CC\oplus\CC$ and define an algebra structure on
this space. Let the element $c\equiv(0,1,0)$ commutes with
everything and the commutator of the element $d\equiv(0,0,1)$ with
the elements $\hat x[s]\equiv(x\otimes s,0,0)$, $x\in\mathfrak a$,
$s\in\mathfrak K$, is given by the formula
$[d,\hx[s]]=\hx[s']$, 
where $s'$ is a derivation of $s$.
Define the Lie bracket between the elements of type $\hx[s]$ requiring the scalar product defined by formulae
\begin{align*}
 \la \hat x[s],\hat y[t]\ra&=\la x\otimes s,y\otimes t\ra, &
 \la c,\hx[s]\ra&=\la d,\hx[s]\ra=0, &
 \la c,d\ra&=1
\end{align*}
to be invariant. It gives the formula
\begin{align}
 [\hat x[s],\hat y[t]]=([x\otimes s,y\otimes t]_0,0,0)+ c\cdot B(x_1\otimes s_1,x_2\otimes s_2), \label{XXc}
\end{align}
where $[\cdot,\cdot]_0$ is the Lie bracket in the algebra
$\frg=\mathfrak a\otimes\mathfrak K$ and $B(\cdot,\cdot)$ is a
standard 1-cocycle:
$B\big(x\otimes s,y\otimes t\big)=(x,y)\la s'(z),t(z)\ra_z$. 
The expression $\hat x[s]$ depends linearly on $s\in\fK$ and,
therefore, can be regarded as an action of operator
$\hx\colon\fK\to\hg$. The commutation relations for the algebra
$\hg$ in terms of currents $x(u)$ corresponding to these operators
can be written in the standard form: $[c,x(u)]=[c,d]=0$ and
\begin{align}
 [x_1(u),x_2(v)]&=x_3(u)\delta(u,v)-c\cdot(x_1,x_2){d\delta(u,v)}/{du}, & [d,x(u)]&=-d x(u)/du, \label{xxtc}
\end{align}
where $x_1,x_2\in\mathfrak a$, $x_3=[x_1,x_2]_{\mathfrak a}$.

\bb{Half-currents.} To describe different bialgebra structures in
the current algebras we have to decompose  the currents in these
algebras into difference of the currents which have good
analytical properties in certain domains: $x(u)=x^+(u)-x^-(u)$.
The $\frg$-valued distributions $x^+(u)$, $x^-(u)$ are called {\em
half-currents}. To perform such a decomposition we will use
so-called Green distributions~\cite{ER97}. Let
$\Omega^+,\Omega^-\subset\Sigma\times\Sigma$ be two domain
separated by a hypersurface $\bar\Delta\subset\Sigma\times\Sigma$
which contains the diagonal:
$\Delta=\{(u,u)\mid u\in\Sigma\}\subset\bar\Delta$. Let there exist distributions
$G^+(u,z)$ and $G^-(u,z)$ regular in $\Omega^+$ and $\Omega^-$
respectively such that $ \delta(u,z)=G^+(u,z)-G^-(u,z)$. To define
half-currents corresponding to these Green distributions we
decompose them as $G^+(u,z)=\sum_i\alpha^+_i(u)\beta^+_i(z)$ and
$G^-(u,z)=\sum_i\alpha^-_i(u)\beta^-_i(z)$. Then the half-currents
are defined as $x^+(u)=\sum_i \alpha^+_i(u)\hx[\beta^+_i]$ and
$x^-(u)=\sum_i \alpha^-_i(u)\hx[\beta^-_i]$. This definition does
not depend on a choice of decompositions of the Green
distributions. The half-currents are currents corresponding to the
operators $\hx^\pm=\pm\,\hx\cdot P^\pm$, where $P^\pm[s](z)=\pm\la
G^\pm(u,z),s(u)\ra$, $s\in\fK$. One can express the half-currents
through the current $x(u)$, which we shall call {\em a total
current} in contrast with the half ones:
\begin{align}\label{hc_tc}
 x^+(u)=\la G^+(u,z)x(z)\ra_z,\qquad x^-(u)=\la G^-(u,z)x(z)\ra_z.
\end{align}
Here $\la a(z) \ra_z\equiv \la a(z),1 \ra_z$.

\bb{Two elliptic classical current algebras.} In this paper we
will consider the case when $\Sigma$ is a covering of an elliptic
curve and Green distributions are regularization of certain
quasi-doubly periodic meromorphic functions. We will call the
corresponding centrally extended algebras of currents by {\em
elliptic classical current algebras}. The main aim of this paper
is to show the following facts:
\begin{itemize}
\item There are two essentially different choices of the test function algebras $\fK$ in this case
corresponding to the different covering $\Sigma$.
\item The same quasi-doubly periodic meromorphic functions regularized with respect to the different
 test function algebras define the different quasi-Lie bialgebra structures and, therefore,
 the different classical elliptic current algebras.
\item The internal structure of these two elliptic algebras is essentially different
in spite of a similarity in the commutation relations between
their half-currents.
\end{itemize}

The first choice corresponds to $\fK=\lfK_0$, where $\lfK_0$
consists of complex-valued one-variable functions defined in a
vicinity of origin (see details in Appendix~\ref{apA}) equipped
with the scalar product~\r{lfK_sp}. These functions can be
extended up to meromorphic functions on the covering
$\Sigma=\mathbb C$. The regularization domain $\Omega^+$,
$\Omega^-$ for Green distributions in this case consist of the
pairs $(u,z)$ such that $\min(1,|\tau|)>|u|>|z|>0$ and
$0<|u|<|z|<\min(1,|\tau|)$ respectively, where $\tau$ is an
elliptic module, and $\bar\Delta=\{(u,z)\mid |u|=|z|\}$.

The second choice corresponds to $\fK=K=K(\Cyl)$. The algebra $K$ consists of entire periodic functions
$s(u)=s(u+1)$ on $\CC$ decaying exponentially at $\Im u\to\pm\infty$ equipped with an invariant scalar
product~\r{spJ}. This functions can be regarded as functions on cylinder $\Sigma=\Cyl$
(see Appendix~\ref{apA}). The regularization domains $\Omega^+$, $\Omega^-$ for Green distributions
consist of the pairs $(u,z)$ such that $-\Im\tau<\Im(u-z)<0$ and $0<\Im (u-z)<\Im\tau$ respectively
and $\bar\Delta=\{(u,z)\mid \Im u=\Im z\}$.

\bb{Integration contour.} The geometric roots of the difference
between these two choices can be explained as follows. These
choices of test functions on different coverings $\Sigma$ of
elliptic curve correspond to the homotopically different contours
on the elliptic curve. Each test function can be considered as an
analytical continuation of a function from this contour -- a real
manifold -- to the corresponding covering. This covering should be
chosen as a most homotopically simple to cover this contour and to
obtain a bigger source of test functions. In the first case, this
contour is a homotopically trivial and coincides with a small
contour around fixed point on the torus. We can always choose a
local coordinate $u$ such that $u=0$ in this point. This explains
the notation $\lfK_0$. This contour corresponds to the covering
$\Sigma=\CC$ and it enters in the pairing~\r{lfK_sp}. In the
second case, it goes along a cycle and it can not be represented
as a closed contour on $\CC$. Hence the most simple covering in
this case is a cylinder $\Sigma=\Cyl$ and the contour is that one
in the pairing~\r{spJ}. This leads to essentially different
properties of the current elliptic algebras based on the test
function algebras $\fK=\lfK_0$ and $\fK=K(\Cyl)$.

\bb{Restriction to the $\slt$ case.} To make these differences more transparent we shall consider
only the simplest case of Lie algebra $\mathfrak a=\slt$ defined as a three-dimensional complex
Lie algebra with commutation relations $[h,e]=2e$, $[h,f]=-2f$ and $[e,f]=h$. We denote the
constructed current algebra $\hg$ for the case $\fK=\lfK_0$ as $\cEF$ and for $\fK=K=K(\Cyl)$ as
$\utau$. These current algebras may be identified with classical limits of the quantized currents
algebra $E_{\tau,\eta}(\slt)$ of  \cite{EF98} and $U_{p,q}(\widehat{\mathfrak{sl}}_2)$
of \cite{JKOS} respectively. The Green
distributions appear in the algebras  $\cEF$ and $\utau$ as a regularization of the same meromorphic
quasi-doubly periodic functions but in different spaces: $(\lfK_0\otimes\lfK_0)'$ and $(K\otimes K)'$
respectively. Primes mean the extension to the space of the distributions. We call them {\em elliptic
Green distributions}. We define the algebras $\cEF$ and $\utau$ to be {\itshape a priori} different,
because the main component of our construction, elliptic Green distributions are {\itshape a priori}
different being understood as distributions of different types: acting in the algebras $\lfK_0$ and $K$
respectively. It means, in particular, that their quantum analogs, the algebras
$E_{\tau,\eta}(\mathfrak{sl}_2)$ and $U_{p,q}(\widehat{\mathfrak{sl}}_2)$ are different.

\section{Half-currents and co-structures}
\label{sec3}

We start with a suitable definition of theta-functions and a
conventional choice of standard bases. This choice is motivated
and corresponds to definitions and notations of ~\cite{EPR}.

\bb{Theta-function.} Let $\tau\in\CC$, $\Im\tau>0$ be a module of
the elliptic curve $\mathbb C/\Gamma$, where $\Gamma=\mathbb
Z+\tau\mathbb Z$ is a period lattice. The odd theta function
$\theta(u)=-\theta(-u)$ is defined as a holomorphic function on
$\CC$ with the properties
\begin{align}
 \theta(u+1)&=-\theta(u), &\theta(u+\tau)&=-e^{-2\pi i u-\pi i\tau}\theta(u), &\theta'(0)&=1.
\end{align}

\subsection{Elliptic Green distributions on $\lfK_0$}
\label{subsec31}

\bb{Dual bases.} Fix a complex number $\lambda$.  Consider the following bases in $\lfK_0$ ($n\geq0$):
$\epsilon_{n;\lambda}(u)=(-u)^n$, $\epsilon^{-n-1;\lambda}(u)=u^n$,
\begin{align}
 \epsilon^{n;\lambda}(u)=\frac1{n!}\left(\frac{\theta(u+\lambda)}
 {\theta(u)\theta(\lambda)}\right)^{(n)}\,,\quad
 \epsilon_{-n-1;\lambda}(u)=\frac{(-1)^n}{n!}\left(\frac{\theta(u-\lambda)}
 {\theta(u)\theta(-\lambda)}\right)^{(n)}\,,  \nn
\end{align}
for $\lambda\not\in\Gamma$ and the bases $\epsilon_{n;0}(u)=(-u)^n$, $\epsilon^{-n-1;0}(u)=u^n$,
\begin{align}
 \epsilon^{n;0}(u)=\frac1{n!}\left(\frac{\theta'(u)}{\theta(u)}\right)^{(n)}\,,\quad
 \epsilon_{-n-1;0}(u)=\frac{(-1)^n}{n!}\left(\frac{\theta'(u)}
 {\theta(u)}\right)^{(n)}\,, \nn
\end{align}
for $\lambda=0$. Here $(\cdot)^{(n)}$ means $n$-times derivation. These bases are dual:
$\la\epsilon^{n;\lambda}(u),\epsilon_{m;\lambda}(u)\ra=\delta^n_m$ and
$\la\epsilon^{n;0}(u),\epsilon_{m;0}(u)\ra=\delta^n_m$
with respect to the scalar product~\r{lfK_sp}
which means
\begin{align}
 &\sum_{n\in\mathbb Z}\epsilon^{n;\lambda}(u)\epsilon_{n;\lambda}(z)=\delta(u,z),
 & &\sum_{n\in\mathbb Z}\epsilon^{n;0}(u),\epsilon_{n;0}(z)=\delta(u,z). \label{dbs}
\end{align}

\bb{Green distributions for $\lfK_0$ and the addition theorems.}
Here we follow the ideas of ~\cite{ER97} and ~\cite{EPR}. We
define the following distribution
\begin{align}
 \G+(u,z)&=\sum_{n\ge0}\epsilon^{n;\lambda}(u)\epsilon_{n;\lambda}(z)\,,\quad
 \G-(u,z)=-\sum_{n<0}\epsilon^{n;\lambda}(u)\epsilon_{n;\lambda}(z)\,,  \label{Glpmdef} \\
 G(u,z)&=\sum_{n\ge0}\epsilon^{n;0}(u)\epsilon_{n;0}(z)=\sum_{n<0}\epsilon^{n;0}(z)\epsilon_{n;0}(u)\,.
 \label{Gdef} \end{align}
One can check that these series converge in sense of distributions and, therefore, define
continuous functionals on $\lfK_0$ called Green distributions. Their action on a test function $s(u)$ reads
\begin{align}
\la\G\pm(u,z),s(u)\ra_u&=
  \oint\limits_{|u|>|z|\atop |u|<|z|}\frac{du}{2\pi i}\frac{\theta(u-z+\lambda)}{\theta(u-z)
  \theta(\lambda)}s(u)\,,
  \label{Glpmact} \\
 \la G(u,z),s(u)\ra_u&=
  \oint\limits_{|u|>|z|}\frac{du}{2\pi i}\frac{\theta'(u-z)}{\theta(u-z)}s(u)\,, \label{Gact}
\end{align}
where integrations are taken over circles around zero which are
small enough such that the corresponding inequality takes place.
The Green distributions are examples of the `shifted'
distributions (see Appendix~\ref{apA}). The formulae~\r{dbs} and
definitions~\r{Glpmdef}, \r{Gdef} imply that
\begin{align}
 \G+(u,z)-\G-(u,z)=\delta(u,z)\,,\quad
 G(u,z)+G(z,u)=\delta(u,z)\,. \label{GGdelta}
\end{align}
The last formulae can be also obtained from~\r{Glpmact},
\r{Gact} taking into account that the function $s(u)$ has poles only in the points $u=0$.
As it is seen from~\r{Glpmact}, the oddness of function $\theta(u)$
leads to the following connection between the $\lambda$-depending Green distributions:
$\G+(u,z)=-G_{-\lambda}^-(z,u)$. 

\begin{prop} \label{prop1}
The {\em semi-direct products} (see Appendix~\ref{apA}) of Green
distributions are related by the following addition formulae
\begin{align}
\begin{split}
 \G+(u,z)\G-(z,v)&=\G+(u,v)G(u,z)-\G+(u,v)G(v,z)-\frac{\partial}{\partial\lambda}\G+(u,v)\,,
\end{split} \label{llpp} \\
\begin{split}
 \G+(u,z)\G+(z,v)&=\G+(u,v)G(u,z)+\G+(u,v)G(z,v)-\frac{\partial}{\partial\lambda}\G+(u,v)\,,
\end{split} \label{llpm} \\
\begin{split}
 \G-(u,z)\G-(z,v)&=-\G-(u,v)G(z,u)-\G-(u,v)G(v,z)-\frac{\partial}{\partial\lambda}\G+(u,v)\,,
\end{split} \label{llmp} \\
\begin{split}
 \G-(u,z)\G+(z,v)&=-\G+(u,v)G(z,u)+\G+(u,v)G(z,v)-\frac{\partial}{\partial\lambda}\G+(u,v)\,.
\end{split} \label{llmm}
\end{align}
\end{prop}

\noindent{\it Proof.}\
 The actions of both hand sides of~\r{llpp}, for example, can be reduced to the integration
 over the same contours with some kernels. One can check the equality of these kernels using the
 degenerated Fay's identity~\cite{Fay}
\begin{align}
\frac{\theta(u-z+\lambda)}{\theta(u-z)\theta(\lambda)} \frac{\theta(z+\lambda)}{\theta(z)\theta(\lambda)}=
 \frac{\theta(u+\lambda)}{\theta(u)\theta(\lambda)} \frac{\theta'(u-z)}{\theta(u-z)}
 +\frac{\theta(u+\lambda)}{\theta(u)\theta(\lambda)} \frac{\theta'(z)}{\theta(z)}
 -\frac{\partial}{\partial\lambda} \frac{\theta(u+\lambda)}{\theta(u)\theta(\lambda)}. \label{Fayid}
\end{align}
The other formulae can be proved in the same way
if one takes into account~\r{a_delta} and~\r{GGdelta}. \qed

\bb{Projections.} Let us notice that the vectors
$\epsilon_{n;\lambda}(u)$ and $\epsilon^{-n-1;\lambda}(u)$ span
two complementary subspaces  of $\lfK_0$. The
formulae~\eqref{Glpmdef} mean that the distributions $\G+(u,z)$
and $\G-(u,z)$ define orthogonal projections $P_\lambda^+$ and
$P_\lambda^-$ onto these subspaces. They act as
$P_\lambda^+[s](z)=\la\G+(u,z),s(u)\ra_u$ and   $P_\lambda^-[s](z)=-\la\G-(u,z),s(u)\ra_u$. 
Similarly, the operators
$P^+[s](z)=\la G(u,z),s(u)\ra_u$ and   $P^-[s](z)=\la G(z,u),s(u)\ra_u$ 
are projections onto the lagrangian (involutive) subspaces spanned
by vectors $\epsilon_{n;0}(u)$ and $\epsilon_{-n-1;0}(u)$,
respectively. The fact that the corresponding spaces are
complementary to each other is encoded in the
formulae~\eqref{GGdelta}, which can be rewritten as
$P_\lambda^++P_\lambda^-=\id$, $ P^++P^-=\id$. The idempotent
properties and orthogonality of these projection
\begin{align*}
P_\lambda^\pm\cdot P_\lambda^\pm&=P_\lambda^\pm, &
P^\pm\cdot P^\pm&=P^\pm, &  P_\lambda^+\cdot P_\lambda^-&=P_\lambda^-\cdot P_\lambda^+\hm=0, &
P^+\cdot P^-&=P^-\cdot P^+&=0
\end{align*}
are encoded in the formulae
\begin{align}
 \la\G+(u,z)\G+(z,v)\ra_z&=\G+(u,v)\,, & \la\G+(u,z)\G-(z,v)\ra_z&=0\,, \label{GpGpconv}\\
 \la\G-(u,z)\G-(z,v)\ra_z&=-\G-(u,v)\,, & \la\G-(u,z)\G+(z,v)\ra_z&=0\,, \label{GmGmconv} \\
 \la G(u,z)G(z,v)\ra_z&=G(u,v)\,, & \la G(u,z)G(v,z)\ra_z&=0\,, \label{GGconv}
\end{align}
which immediately follow from~\eqref{Glpmdef}, \eqref{Gdef} and also can be obtained from the
relations~\eqref{llpp} -- \eqref{llmm} if one takes into account $\la G(u,z)\ra_z=0$, $\la G(z,u)\ra_z=1$.

\subsection{Elliptic Green distributions on $K$}
\label{subsec32}

\bb{Green distributions and dual bases for $K$.} The analogs of
the Green distributions $\G+(u,z)$, $\G-(u,z)$ are defined in this case
by the following action on the space~$K$
\begin{align}
 \la\Gg\pm(u-z),s(u)\ra_u&=\int\limits_{-\Im\tau<\Im(u-z)<0\atop 0<\Im(u-z)<\Im\tau}
 \frac{du}{2\pi i} \frac{\theta(u-z+\lambda)}{\theta(u-z)\theta(\lambda)}s(u), \label{GlpmdefJ} \\
 \la{\cal G}(u-z),s(u)\ra_u
   &=\int\limits_{-\Im\tau<\Im(u-z)<0}\frac{du}{2\pi i}\frac{\theta'(u-z)}{\theta(u-z)}s(u). \label{GdefJ}
\end{align}
where we integrate over line segments of unit length (cycles of cylinder) such that
the corresponding inequality takes place. The role of dual bases in the algebra $K$ is played by
$\{j_n(u)=e^{2\pi inu}\}_{n\in\mathbb Z}$ and
$\{j^n(u)=2\pi i e^{-2\pi inu}\}_{n\in\mathbb Z}$ (see Appendix~\ref{apA}), a decomposition to
these bases is the usual Fourier expansion. The Fourier expansions for the Green distributions
are~\footnote{Fourier expansions presented in this subsection are obtained considering integration
around boundary of fundamental domain (see~\cite{WW}).}
\begin{gather}
 \Gg\pm(u-z)=\pm2\pi i\sum_{n\in\mathbb Z}\frac{e^{-2\pi in(u-z)}}{1-e^{\pm2\pi i(n\tau-\lambda)}},
 \label{Ggpmdec} \\
{\cal G}(u-z)=\pi i+2\pi i\sum_{n\ne0}\frac{e^{-2\pi in(u-z)}}{1-e^{2\pi in\tau}}. \label{GexpanJ}
\end{gather}
These expansions are in according with formulae
\begin{align}
 \Gg+(u-z)-\Gg-(u-z)&=\delta(u-z), \label{GpGmdeltaJ} \\
  {\cal G}(u-z)+{\cal G}(z-u)&=\delta(u-z), \label{GGdeltaJ}
\end{align}
where $\delta(u-z)$ is a delta-function on $K$, given by the
expansion~\r{dffeJ}.

\bb{Addition formulae.} Now we obtain some properties of these
Green distributions and compare them with the properties of their
analogs $\G+(u,z)$, $\G-(u,z)$, $G(u,z)$ described in
subsection~\ref{subsec31}. In particular, we shall see that some
properties are essentially different. Let us start with the
properties of Green distribution which are similar to the case of
algebra $\lfK_0$. They satisfy the same addition formulae that was
described in~\ref{subsec31}:
\begin{prop} \label{prop2}
 The semi-direct product of Green distributions for algebra $K$ is related by the
 formulae~\r{llpp}--\r{llmm} with the distributions $\Gg\pm(u-z)$,
 ${\cal G}(u-z)$ instead of $\G\pm(u-z)$, $G(u-z)$ respectively.
\end{prop}

\noindent{\it Proof.}\  The kernels of these distributions are the same and therefore
the addition formula in this case is also based on the Fay's identity~\r{Fayid}. \qed \\

\bb{Analogs of projections.} The Green distributions define the following operators on~$K$:
\begin{align*}
 {\cal P}_\lambda^+[s](z)&=\la\Gg+(u-z),s(u)\ra_u, & {\cal P}_\lambda^-[s](z)&=\la\Gg-(u-z),s(u)\ra_u, \\
  {\cal P}^+[s](z)&=\la {\cal G}(u-z),s(u)\ra_u, & {\cal P}^-[s](z)&=\la {\cal G}(z-u),s(u)\ra_u
\end{align*}
which are similar to their analogs $P_\lambda^\pm$, $P^\pm$ and also satisfy ${\cal P}_\lambda^++{\cal P}_\lambda^-=\id$, ${\cal P}^++{\cal P}^-=\id$ (due to~\r{GpGmdeltaJ}), but they are not projections (idempotent). This fact is reflected in the following relations, which are consequence of the formulae~\r{llpp}--\r{llmm} and
$\la{\cal G}(u-z)\ra_z=\la{\cal G}(z-u)\ra_z=1/2$:
\begin{align}
 \la\Gg+(u-z)\Gg+(z-v)\ra_z&=\Gg+(u-v)-\frac1{2\pi i}\frac{\partial}{\partial\lambda}\Gg+(u-v)\,,
 \label{GpGpconvJ} \\
 \la\Gg+(u-z)\Gg-(z-v)\ra_z&=-\frac1{2\pi i}\frac{\partial}{\partial\lambda}\Gg+(u-v)\,,
 \label{GpGmconvJ} \\
 \la\Gg-(u-z)\Gg+(z-v)\ra_z&=-\frac1{2\pi i}\frac{\partial}{\partial\lambda}\Gg+(u-v)\,,
 \nn\\ 
 \la\Gg-(u-z)\Gg-(z-v)\ra_z&=-\Gg-(u-v)-\frac1{2\pi i}\frac{\partial}{\partial\lambda}\Gg+(u-v)\,,
 \label{GmGmconvJ}
\end{align}
where $\frac{\partial}{\partial\lambda}\Gg+(u-v)=\frac{\partial}{\partial\lambda}\Gg-(u-v)$ and
\begin{align}
\la{\cal G}(u-z){\cal G}(z-v)\ra_z=&{\cal G}(u-v)
 -\frac1{4\pi i}\gamma(u-v)\,, \label{GGppconvJ} \\
\la{\cal G}(u-z){\cal G}(v-z)\ra_z=&\la{\cal G}(z-u){\cal G}(z-v)\ra_z=
 \frac1{4\pi i}\gamma(u-v)\,. \label{GGpmconvJ}
\end{align}
$\gamma(u-z)$ is a distribution which has the following action and expansion
\begin{align}
 \la&\gamma(u-z),s(u)\ra=-\frac{\theta'''(0)+4\pi^2}3+\int\limits_{-\frac12}^{+\frac12}\frac{du}{2\pi i}
 \frac{\theta''(u-z)}{\theta(u-z)}s(u)\,, \nn\\
&\gamma(u-z)=-2\pi^2+8\pi^2\sum_{n\ne0}\frac{e^{-2\pi in(u-z)+2\pi in\tau}}{(1-e^{2\pi in\tau})^2}\,.
\nn
\end{align}

\bb{Comparison of the Green distributions.} Contrary
to~\r{GpGpconv}, \r{GmGmconv} the
formulae~\eqref{GpGpconvJ}--\eqref{GGpmconvJ} contain some
additional terms in the right hand sides obstructed the operators
${\cal P}_\lambda^\pm$, ${\cal P}^\pm$ to be projections. They do
not decompose the space $K(\Cyl)$ in a direct sum of subspaces
as it would be in the case of projections $P_\lambda^\pm$, $P^\pm$
acting on $\lfK_0$. Moreover, as one can see from the Fourier
expansions~\r{Ggpmdec}, \r{GexpanJ} of Green distributions the
images of the operators coincide with whole algebra $K$:
${\cal P}_\lambda^\pm\big(K(\Cyl)\big)=K(\Cyl)$, ${\cal
P}^\pm\big(K(\Cyl)\big)=K(\Cyl)$. As we shall see this fact has a
deep consequence for the half-currents of the corresponding Lie
algebra $\utau$. As soon as we are aware that the positive
operators ${\cal P}_\lambda^+$, ${\cal P}^+$ as well as negative
ones ${\cal P}_\lambda^-$, ${\cal P}^-$
 transform the algebra $K$ to itself, we can surmise that
they can be related to each other. This is actually true. From formulae~\r{Ggpmdec}, \r{GexpanJ} we
conclude that
\begin{align}
\Gg+(u-z-\tau)&=e^{2\pi i\lambda}\Gg-(u-z),  &{\cal G}(u-z-\tau)&=2\pi i-{\cal G}(z-u). \label{shift_tau}
\end{align}
In terms of operator's composition these properties look as
\begin{align}
 {\cal T}_\tau\circ{\cal P}_\lambda^+&={\cal P}_\lambda^+\circ{\cal T}_\tau=-e^{2\pi i\lambda}
{\cal P}_\lambda^-, &
{\cal T}_\tau\circ{\cal P}^+&={\cal P}^+\circ{\cal T}_\tau=2\pi i{\cal I}-{\cal P}^-,
\end{align}
where ${\cal T}_t$ is a shift operator: ${\cal T}_t[s](z)=s(z+t)$,
and ${\cal I}$ is an integration operator: ${\cal
I}[s](z)=\int_{-\frac12}^{\frac12}\frac{du}{2\pi i}s(u)$. This
property is no longer true for the case of
 Green distributions from section~\ref{subsec31}.

\subsection{Elliptic half-currents}
\label{subsec33}

\bb{Tensor subscripts.} First introduce the following
notation. Let $\bU={\cal U}(\mathfrak{g})$ be a universal
enveloping algebra of the considering Lie algebra $\mathfrak{g}$
and $V$ be a $\bU$-module. For an element
$t=\sum_{k}a_1^k\otimes\ldots\otimes a_n^{k}\otimes
u_1^k\otimes\ldots\otimes u_m^{k} \in\End V^{\otimes
n}\otimes\bU^{\otimes m}$, where $n,m\ge0$,
$a_1^k,\ldots,a_n^k\in\End V$, $u_1^k,\ldots,u_m^k\in\bU$ we shall
use the following notation for an element of $\End V^{\otimes
N}\otimes\bU^{\otimes M}$, $N\ge n, M\ge m$,
\begin{align*}
 t_{i_1,\ldots,i_n,{\bf j_1},\ldots,{\bf j_m}}=&
  \sum_k \id_V\otimes\cdots\otimes\id_V\otimes a_1^k\otimes\id_V\otimes\cdots\otimes\id_V\otimes
  a_n^k\otimes\id_V\otimes\cdots\otimes\id_V\otimes \\
  &\otimes1\otimes\ldots\otimes1\otimes
  u_1^k\otimes1\otimes\ldots\ldots\otimes u_m^{k}
  \otimes1\otimes\ldots\otimes1,
\end{align*}
where $a_s^k$ stays in the $i_s$-th position in the tensor product and $u_s^k$ stays in the $j_s$-th position. \\

\bb{Half-currents.} The total currents $h(u)$, $e(u)$ and $f(u)$ of the algebra $\cEF$ can be divided into
half-currents using the Green distributions $G(u,z)$, $-G(z,u)$ for $h(u)$; $\G+(u,z)$,
$\G-(u,z)$ for $e(u)$; and $G^+_{-\lambda}(u,z)=-\G-(z,u)$, $G^-_{-\lambda}(u,z)=-\G+(z,u)$ for $f(u)$.
The relations of type~\r{hc_tc}, then, looks as
\begin{align}
 h^+(u)&=\la G(u,v)h(v)\ra_v, & h^-(u)&=-\la G(v,u)h(v)\ra_v,         \label{hpmht} \\
 \e^+(u)&=\la \G+(u,v)e(v)\ra_v, & \e^-(u)&=\la \G-(u,v)e(v)\ra_v,     \label{epmht} \\
 \f^+(u)&=\la G_{-\lambda}^+(u,v)f(v)\ra_v, &  \f^-(u)&=\la G_{-\lambda}^-(u,v)f(v)\ra_v,  \label{fpmht}
\end{align}
so that
  $h(u)=h^+(u)-h^-(u)$, $e(u)=\e^+(u)-\e^-(u)$, $f(u)=\f^+(u)-\f^-(u)$.

\bb{$rLL$-relations for $\cEF$.} The commutation relation between
the half-currents can be written in a matrix form. Let us
introduce the matrices of {\em $L$-operators:}
\begin{gather}
\L_\lambda^\pm(u)=
 \begin{pmatrix}
  \frac12h^\pm(u) & \f^\pm(u) \\
  \e^\pm(u)       & -\frac12h^\pm(u)
 \end{pmatrix}, \\
\end{gather}
as well as the $r$-matrices:
\begin{gather} r_\lambda^+(u,v)=
 \begin{pmatrix}
   \frac12G(u,v) & 0              & 0              & 0              \\
  0              & -\frac12G(u,v) & G^+_{-\lambda}(u,v)      & 0              \\
  0              & \G+(u,v)       & -\frac12G(u,v) & 0              \\
  0              & 0              & 0              & \frac12G(u,v)
 \end{pmatrix}.
\end{gather}

\begin{prop} \label{prop3}
 The commutation relations of the algebra $\cEF$ in terms of half-currents
 can be written in the form:
\begin{align}
   [d,\L_{\lambda}^\pm(u)]=-\frac{\partial}{\partial u}\L_{\lambda}^\pm(u), \label{dL}
\end{align}
\begin{multline} \label{rLpmLpm}
 [\L_{\lambda,1}^\pm(u),\L_{\lambda,2}^\pm(v)]=[\L_{\lambda,1}^\pm(u)+\L_{\lambda,2}^\pm
 (v),r_\lambda^+(u-v)]+ \\
 +H_1\frac{\partial}{\partial\lambda}\L_{\lambda,2}^\pm(v)-H_2\frac{\partial}{\partial\lambda}
 \L_{\lambda,1}^\pm(u)
 +h\frac{\partial}{\partial\lambda}r_\lambda^+(u-v),
\end{multline}
\begin{multline} \label{rLpLm}
 [\L_{\lambda,1}^+(u),\L_{\lambda,2}^-(v)]=[\L_{\lambda,1}^+(u)+\L_{\lambda,2}^-(v),r_\lambda^+(u-v)]+ \\
 +H_1\frac{\partial}{\partial\lambda}\L_{\lambda,2}^-(v)-H_2\frac{\partial}{\partial
 \lambda}\L_{\lambda,1}^+(u)
 +h\frac{\partial}{\partial\lambda}r_\lambda^+(u-v)+c\cdot\frac{\partial}{\partial u}r_\lambda^+(u-v),
\end{multline}
where
$H=
 \begin{pmatrix}
  1 & 0 \\
  0  & -1
 \end{pmatrix}$ and  $h=\hat h[\epsilon_{0;0}]$.
The $L$-operators satisfy an important relation
\begin{align}
 [H+h,\L^\pm_\lambda(u)]=0\,. \label{HhL}
\end{align}
\end{prop}

\noindent{\it Proof.}\  Using the formulae~\eqref{GpGpconv} -- \eqref{GGconv} we calculate the
scalar products on the half-currents:
  $\la \L_{\lambda,1}^\pm(u),\L_{\lambda,2}^\pm(v)\ra=0$,
    $\la \L_{\lambda,1}^+(u),\L_{\lambda,2}^-(v)\ra=-r_\lambda^+(u,v)$.
Differentiating these formulae by $u$  we can  obtain the values of the standard co-cycle
on the half-currents:
  $B\big(\L_{\lambda,1}^\pm(u),\L_{\lambda,2}^\pm(v)\big)=0$,
  $B\big(\L_{\lambda,1}^+(u),\L_{\lambda,2}^-(v)\big)=\frac{\partial}{\partial u}r_\lambda^+(u,v)$.
Using the formulae~\r{llpp}--\r{llmm} one can calculate the
brackets $[\cdot,\cdot]_0$ on the half-currents. Representing them
in the matrix form and adding the co-cycle term one can derive the
relations~\r{rLpmLpm}, \eqref{rLpLm}. Using the formulae
$[h,\L_{\lambda}^\pm(v)]=\tr_1\la
H_1[\L_{\lambda,1}^+(u),\L_{\lambda,2}^\pm(v)]\ra_u$, $\tr_1\la
H_1r_{\lambda}^+(u,v)\ra_u=H$,
$\tr_1\la[H_1,\L_{\lambda,1}^+(u)]r_{\lambda}^+(u,v)\ra_u=0$ we
obtain the relation~\eqref{HhL} from~\eqref{rLpmLpm},
\eqref{rLpLm}. \qed

\bb{$rLL$-relations for $\utau$.} Now consider the case of the
algebra $\utau$. The half-currents, $L$-operators
$\LL_\lambda^\pm(u)$ and $r$-matrix $\rr_\lambda^+(u-v)$ are
defined by the same formulae as above with distributions
${G}(u,v)$ and ${G}^\pm_{\lambda}(u,v)$ replaced everywhere by the
distributions ${\cal G}(u,v)$  and ${\cal G}^\pm_{\lambda}(u,v)$.
We have
\begin{prop} \label{prop4}
The commutation relations of the algebra $\cEF$ in terms of
half-currents  can be written in the form:
\begin{multline}
 [\LL_{\lambda,1}^\pm(u),\LL_{\lambda,2}^\pm(v)]=[\LL_{\lambda,1}^\pm(u)+
 \LL_{\lambda,2}^\pm(v),\rr_\lambda^+(u-v)]+ \\
 +H_1\frac{\partial}{\partial\lambda}\LL_{\lambda,2}^\pm(v)-H_2\frac{\partial}
 {\partial\lambda}\LL_{\lambda,1}^\pm(u)
 +h\frac{\partial}{\partial\lambda}\rr_\lambda^+(u-v)-c\cdot\frac{\partial}
 {\partial\tau}\rr_\lambda^+(u-v)\,, \label{rLpmLpmJ}
\end{multline}
\begin{multline}
 [\LL_{\lambda,1}^+(u),\LL_{\lambda,2}^-(v)]=[\LL_{\lambda,1}^+(u)+
 \LL_{\lambda,2}^-(v),\rr_\lambda^+(u-v)]+ \\
 +H_1\frac{\partial}{\partial\lambda}\LL_{\lambda,2}^-(v)-H_2\frac{\partial}
 {\partial\lambda}\LL_{\lambda,1}^+(u)
 +h\frac{\partial}{\partial\lambda}\rr_\lambda^+(u-v)
 +c\cdot\bigg(\frac{\partial}{\partial u}-\frac{\partial}{\partial\tau}\bigg)
 \rr_\lambda^+(u-v)\,, \label{rLpLmJ}
\end{multline}
where $h=\hat h[j_0]$. We also have in this case the relation
\begin{align}
 [H+h,\LL^\pm_\lambda(u)]=0\,. \label{HhLJ}
\end{align}
\end{prop}

\noindent{\it Proof.}\  To express the standard co-cycle on the half currents through
the derivatives of the $r$-matrix we need the following formulae
\begin{align}
 \frac1{2\pi i}\frac{\partial}{\partial u}\frac{\partial}{\partial\lambda}\Gg+(u-v)
     &=\frac{\partial}{\partial\tau}\Gg+(u-v), \nn\\ 
 \frac1{2\pi i}\frac{\partial}{\partial u}\frac{\partial}{\partial\lambda}\Gg-(u-v)
     &=\frac{\partial}{\partial\tau}\Gg-(u-v)=\frac{\partial}{\partial\tau}\Gg+(u-v),
     \nn\\ 
 \frac1{4\pi i}\frac{\partial}{\partial u}\gamma(u-v)&=\frac{\partial}{\partial\tau}{\cal G}(u-v).
 \nn  
\end{align}
Using these formulae we obtain
\begin{align}
  B\big(\LL_{\lambda,1}^\pm(u),\LL_{\lambda,2}^\pm(v)\big)&=-\frac{\partial}{\partial\tau}\rr_\lambda^+(u-v),
    & B\big(\LL_{\lambda,1}^+(u),\LL_{\lambda,2}^-(v)\big)
      &=\bigg(\frac{\partial}{\partial u}-\frac{\partial}{\partial\tau}\bigg)\rr_\lambda^+(u-v).
\end{align}
Using the formulae $[h,\LL_{\lambda}^\pm(v)]=2\tr_1\la
H_1[\LL_{\lambda,1}^+(u),\L_{\lambda,2}^\pm(v)]\ra_u$, $\tr_1\la
H_1r_{\lambda}^+(u,v)\ra_u=H/2$,
$\tr_1\la[H_1,\L_{\lambda,1}^+(u)]r_{\lambda}^+(u,v)\ra_u
=\frac{i}{\pi
}\frac{\partial}{\partial\lambda}\LL_{\lambda}^\pm(v)$ we get the
relation~\r{HhLJ} from~\r{rLpmLpmJ}, \r{rLpLmJ}. \qed

\bb{Peculiarities of half-currents for $\utau$. } To conclude this subsection we discuss the
implication of  the properties of the Green distributions described in the end of the previous section
to the Lie algebra $\utau$. The fact that the images of the operators ${\cal P}_\lambda^\pm$,
${\cal P}^\pm$ coincide with all the space $K$ means that the
commutation relations between the positive (or negative) half-currents are sufficient to describe
all the Lie algebra $\utau$.
This is a consequence of construction of the Lie algebra $\utau$ as the
central extension of the loop algebra $\slt\otimes K$. To obtain all commutation relations given in
Proposition~\ref{prop4}
from relations between only positive (or negative) half-currents one can use, first,
the connection between positive and negative half-currents:
\begin{align*}
 h^+(u-\tau)&=2\pi i h+h^-(u), & e^+(u-\tau)&=e^{2\pi i\lambda}e^-(u), & f^+(u-\tau)&=e^{-2\pi i\lambda}f^-(u),
\end{align*}
which follows from the properties of Green distributions expressed in formulae~\r{shift_tau}.
Second, relations~\r{HhLJ}, which also follow from the relations between only positive (respectively
negative) half-currents, and finally, one needs to use the equality
$\frac{\partial}{\partial\tau}\Gg\pm(u-z-\tau)=
  e^{2\pi i\lambda}(-\frac{\partial}{\partial u}+\frac{\partial}{\partial\tau})\Gg-(u-z)$.
At this point we see the essential difference of the Lie algebra $\utau$ with the Lie algebra $\cEF$.

\subsection{Coalgebra  structures of $\cEF$ and $\utau$}
\label{subsec34}
We describe here the structure of {\em quasi-Lie bialgebras} for our Lie algebras $\cEF$ and $\utau$.
We will start with an explicit expression for universal (dynamical) $r$-matrices for both Lie algebras.

\begin{prop} \label{prop_CDYBE}
The universal $r$-matrix for the Lie algebra $\cEF$ defined as
\begin{align}
 r_\lambda=\frac12\sum_{n\ge0}\hat h[\epsilon^{n;0}]\otimes\hat h[\epsilon_{n;0}]
  +\sum_{n\ge0}\hat f[\epsilon^{n;\lambda}]\otimes\hat e[\epsilon_{n;\lambda}]
  +\sum_{n<0}\hat e[\epsilon_{n;\lambda}]\otimes\hat f[\epsilon^{n;\lambda}]+c\otimes d \label{curm}
\end{align}
satisfies the Classical Dynamical Yang-Baxter Equation (CDYBE)
\begin{align}
  [r_{\lambda,{\bf 12}},r_{\lambda,{\bf 13}}]
 +[r_{\lambda,{\bf 12}},r_{\lambda,{\bf 23}}]
 +[r_{\lambda,{\bf 13}},r_{\lambda,{\bf 23}}]=
   h_{{\bf 1}}\frac{\partial}{\partial\lambda}r_{\lambda,{\bf 23}}
  -h_{{\bf 2}}\frac{\partial}{\partial\lambda}r_{\lambda,{\bf 13}}
  +h_{{\bf 3}}\frac{\partial}{\partial\lambda}r_{\lambda,{\bf 12}}\,. \label{CDYBE}
\end{align}
\end{prop}

Denote by $\Pi_u$  the evaluation representation $\Pi_u\colon\cEF\to\End V_u$,
where $V_u=\mathbb C^2\otimes\lfK_0$ and the subscript $u$ means the argument
of the functions belonging to $\lfK_0$:
\begin{align}
 &\Pi_u\colon \hat h[s]\mapsto s(u)H, & &\Pi_u\colon \hat e[s]\mapsto s(u)E, &
 &\Pi_u\colon \hat f[s]\mapsto s(u)F, \label{ev_repr}
\end{align}
and
$\Pi_u\colon c\mapsto 0$, $\Pi_u\colon d\mapsto \frac{\partial}{\partial u}$, where $ H=\left(
 \begin{smallmatrix}
  1 & 0 \\
  0  & -1
 \end{smallmatrix}\right)$,
 $ E=\left(
 \begin{smallmatrix}
  0 & 1 \\
  0  & 0
 \end{smallmatrix}\right)$,
 $ F=\left(
 \begin{smallmatrix}
  0 & 0 \\
  1  & 0
 \end{smallmatrix}\right)$,
$s\in\lfK_0$. The relations between $L$-operators, $r$-matrix and universal $r$-matrix are given by the formulae
\begin{align}
\L_\lambda^+(u)&=(\Pi_u\otimes\id) r_\lambda, & 
\L_\lambda^-(u)-c\frac{\partial}{\partial u}&=-(\Pi_u\otimes\id) r_{\lambda,\textbf{21}},  \label{Lm_ev_r}
\end{align}
and $r_\lambda^+(u-v)=(\Pi_u\otimes\Pi_v) r_\lambda$.
Taking into account these formulae  and applying $(\Pi_u\otimes\Pi_v\otimes\id)$,
$(\id\otimes\Pi_u\otimes\Pi_v)$, $(\Pi_u\otimes\id\otimes\Pi_v)$ to the
equation~\eqref{CDYBE} we derive the relation~\eqref{rLpmLpm} for the sign `$+$',
the relation~\eqref{rLpmLpm} for the sign `$-$' and the relation~\eqref{rLpLm} respectively.
Applying $(\Pi_u\otimes\id)$ or $(\id\otimes\Pi_u)$ to the identity $[\Delta h,\rr_\lambda]=0$
we derive the relation~\eqref{HhL}.

The co-bracket $\delta\colon\cEF\to\cEF\wedge\cEF$ and an element $\varphi\in\cEF\wedge\cEF\wedge\cEF$
are defined as
$\delta x=[\Delta x,r_\lambda]=[x\otimes1+1\otimes x,r_\lambda]$, for
 $x\in\cEF$ and
\begin{align}
 \varphi&=-[r_{\lambda,\textbf{12}},r_{\lambda,\textbf{13}}]
 -[r_{\lambda,\textbf{12}},r_{\textbf{23}}]-[r_{\lambda,\textbf{13}},r_{\lambda,\textbf{23}}]
 = \nn \\
 &=-h_{\textbf{1}}\frac{\partial}{\partial\lambda}r_{\lambda,\textbf{23}}
  +h_{\textbf{2}}\frac{\partial}{\partial\lambda}r_{\lambda,\textbf{13}}
  -h_{\textbf{3}}\frac{\partial}{\partial\lambda}r_{\lambda,\textbf{12}}.\ , \label{varphi_def}
\end{align}
They equip the Lie algebra $\cEF$ with a structure of a quasi-Lie bialgebra~\cite{D90}.
This fact follows from the equality $r_{\textbf{12}}+r_{\textbf{21}}=\Omega$,
where $\Omega$ is a tensor Casimir element of the algebra $\cEF$.
 To calculate this co-bracket on the half-currents in the matrix
 form we apply $(\Pi_u\otimes\id\otimes\id)$, $(\id\otimes\id\otimes\Pi_u)$
 to the equation~\eqref{CDYBE} and derive
\begin{align*}
 \delta L^+_{\lambda}(u)=-[L^+_{\lambda,\textbf{1}}(u)&,L^+_{\lambda,\textbf{2}}(u)]
  +H\frac{\partial}{\partial\lambda}r_{\lambda}
  -h\wedge\frac{\partial}{\partial\lambda}L^+_{\lambda}(u)\,, \\
 \delta L^-_{\lambda}(u)=-[L^-_{\lambda,\textbf{1}}(u)&,L^-_{\lambda,\textbf{2}}(u)]
  +H\frac{\partial}{\partial\lambda}r_{\lambda}
  -h\wedge\frac{\partial}{\partial\lambda}L^-_{\lambda}(u)
  -c\wedge\frac{\partial}{\partial u}L^-_{\lambda}(u)\,.
\end{align*}
We can see also that
 $\delta h=0$, $\delta c=0$, $\delta d=0$.

\begin{prop} \label{prop_CDYBEJ}
The universal $r$-matrix for the Lie algebra $\utau$ defined by formula
\begin{align}
 \rr_\lambda=\frac14\hat h[j^0]\otimes\hat h[j_0]&+\frac12\sum_{n\ne0}\frac{\hat h[j^n]\otimes\hat
  h[j_n]}{1-e^{2\pi in\tau}}+ \notag \\
 &+\sum_{n\in\mathbb Z}\frac{\hat e[j^n]\otimes\hat f[j_n]}{1-e^{2\pi i(n\tau+\lambda)}}
 +\sum_{n\in\mathbb Z}\frac{\hat f[j^n]\otimes\hat e[j_n]}{1-e^{2\pi i(n\tau-\lambda)}}
 +c\otimes d. \label{urm_J}
\end{align}
satisfies the equation
\begin{multline}
  [\rr_{\lambda,\bf{12}},\rr_{\lambda,\bf{13}}]
 +[\rr_{\lambda,\bf{12}},\rr_{\lambda,\bf{23}}]
 +[\rr_{\lambda,\bf{13}},\rr_{\lambda,\bf{23}}]=\\
  =h_{\bf{1}}\frac{\partial}{\partial\lambda}\rr_{\lambda,\bf{23}}
  -h_{\bf{2}}\frac{\partial}{\partial\lambda}\rr_{\lambda,\bf{13}}
  +h_{\bf{3}}\frac{\partial}{\partial\lambda}\rr_{\lambda,\bf{12}}
  -c_{\bf{1}}\frac{\partial}{\partial\tau}\rr_{\lambda,\bf{23}}
  +c_{\bf{2}}\frac{\partial}{\partial\tau}\rr_{\lambda,\bf{13}}
  -c_{\bf{3}}\frac{\partial}{\partial\tau}\rr_{\lambda,\bf{12}}. \label{CDYBEJ}
\end{multline}
\end{prop}

The relations between the universal matrix $\rr_{\lambda}$ and $L$-operators of the
algebra $\utau$ are the same as for the algebra $\cEF$ with a proper modification of the
evaluation representation
$\Pi_u\colon\utau\to\End {\cal V}_u$, ${\cal V}_u=\mathbb C^2\otimes K$ defined by the same
formulas~\r{ev_repr} as above but for
$s\in K$.

The bialgebra structure of $\utau$ is defined in analogous way as for the algebra $\cEF$
and  can be present in the form
\begin{align}
 &\delta\LL^+_{\lambda}(u)=-[\LL^+_{\lambda,\textbf{1}}(u),\LL^+_{\lambda,\textbf{2}}(u)]
  +H\frac{\partial}{\partial\lambda}\rr_{\lambda}
  -h\wedge\frac{\partial}{\partial\lambda}\LL^+_{\lambda}(u)+c\wedge\frac{\partial}{\partial\tau}
  \LL^+_{\lambda}(u)\,, \\
 &\delta\LL^-_{\lambda}(u)=-[\LL^-_{\lambda,\textbf{1}}(u),\LL^-_{\lambda,\textbf{2}}(u)]
  +H\frac{\partial}{\partial\lambda}\rr_{\lambda}
  -h\wedge\frac{\partial}{\partial\lambda}\LL^-_{\lambda}(u)
  -c\wedge\Big(\frac{\partial}{\partial u}-\frac{\partial}{\partial\tau}\Big)\LL^-_{\lambda}(u)\,.
\end{align}

\section{Degenerated cases}
\label{sec4}
We will describe a behavior of our algebras while one or both periods of the elliptic curve become infinite.
The corresponding `degenerated' Green distributions, $r$-matrix and $L$-operators give us
a classical rational or a classical trigonometric `limit' of corresponding elliptic current algebras.

\subsection{Degenerations of the quasi-Lie bialgebra $\cEF$}
\label{subsec41}

There are two different degenerations denoted {\bf(a)} and {\bf(b)} for $\cEF$.
{\bf(a)} corresponds to the case when both periods are infinite ($\omega\to\infty$,
$\omega'\to\infty$). This is a rational degeneration. In the case {\bf(b)} one of the periods
is infinite ($\omega'\to\infty$) while another ($\omega$) rests finite. This is a case of trigonometric
degeneration. A situation when $\omega\to\infty$ and $\omega'$ is finite,
is equivalent to {\bf(b)} due to the symmetry of integration contour and, therefore,
we do not consider it separately.

\bb{Case (a): $\omega\to\infty$, $\omega'\to\infty$,
($\Im\frac{\omega'}{\omega}>0$)}. In order to turn to the lattice
of periods $\Gamma=\mathbb Z \omega+\mathbb Z \omega'$, with
$\frac{\omega'}{\omega}=\tau$, we need to re-scale the variables $u$, $v$ and the dynamical parameter $\lambda$
like $u\to\frac{u}{\omega}$. Let us introduce the following
notations for rational Green distributions
\begin{align*}
\la\ph\pm(u,z),s(u)\ra_u=\oint\limits_{\substack{|u|>|z| \\|u|<|z|}}\frac{du}{2\pi i}\frac1{u-z}s(u).
\end{align*}
These distributions are degenerations of elliptic Green
distributions:
\begin{align*}
 \frac1{\omega}G\big(\frac{u}{\omega},\frac{z}{\omega}\big)&\to\ph+(u,z), &
 \frac1{\omega}G_{\frac\lambda{\omega}}^\pm\big(\frac{u}{\omega},\frac{z}{\omega}\big)&\to
 \frac1{\lambda}+\ph\pm(u,z),
\end{align*}
and the $r$-matrix tends to
\begin{align*}
r_\lambda^{(a)+}(u,v)&=\lim_{\omega,\omega'\to\infty}\frac1{\omega}r_{\frac\lambda{\omega}}^+
\Big(\frac{u}{\omega},\frac{v}{\omega}\Big)= \nn \\
 &=\begin{pmatrix}
   \frac12\ph+(u,v) & 0              & 0              & 0              \\
  0              & -\frac12\ph+(u,v) & -\frac1{\lambda}+\ph+(u,v)      & 0              \\
  0              & \frac1{\lambda}+\ph+(u,v)       & -\frac12\ph+(u,v) & 0              \\
  0              & 0              & 0              & \frac12\ph+(u,v)
 \end{pmatrix}. \nn
\end{align*}
Actually the quasi-Lie bialgebras obtained as rational
degenerations of $\cEF$ for different $\lambda$ is related to each
other by very simple twist. Therefore we shall consider only one
value of the parameter $\lambda$, namely we shall consider the
limited value $\lambda\to\infty$. The $r$-matrix and $L$-operators
looks then as follows
\begin{align*}
 \L^{(a)\pm}(u)&=\lim\limits_{\lambda\to\infty}\lim\limits_{\omega,\omega'\to\infty}\dfrac1{\omega}
   \L_{\frac{\lambda}{\omega}}^\pm\Big(\dfrac{u}{\omega}\Big)=
 \begin{pmatrix}
  \frac12h^{(a)\pm}(u) & f^{(a)\pm}(u) \\
  e^{(a)\pm}(u)       & -\frac12h^{(a)\pm}(u)
 \end{pmatrix},
 \end{align*}
\begin{align}
r^{(a)+}(u,v)&=\lim_{\lambda\to\infty}\frac1{\omega}r_\lambda^{(a)+}(u,v)= \nn \\
 &=\begin{pmatrix}
   \frac12\ph+(u,v) & 0              & 0                 & 0              \\
  0              & -\frac12\ph+(u,v) & \ph+(u,v)         & 0              \\
  0              & \ph+(u,v)         & -\frac12\ph+(u,v) & 0              \\
  0              & 0                 & 0                 & \frac12\ph+(u,v)
 \end{pmatrix}. \label{ra}
\end{align}

Substituting $u\to\frac{u}{\omega}$, $v\to\frac{v}{\omega}$, $\lambda\to\frac{\lambda}{\omega}$
to~\eqref{rLpmLpm} and \eqref{rLpLm} multiplying it by $\frac1{\omega^2}$ and passing to
the limits we obtain
\begin{align}
 [\L_1^{(a),\pm}(u),\L_2^{(a)\pm}(v)]&=[\L_1^{(a)\pm}(u)+\L_2^{(a)\pm}(v),r^{(a)+}(u,v)]. \label{rLpmLpma} \\
 [\L_1^{(a)+}(u),\L_2^{(a)-}(v)]&=[\L_1^{(a)+}(u)+\L_2^{(a)-}(v),r^{(a)+}(u,v)]
   +c\cdot\frac{\partial}{\partial u}r^{(a)+}(u,v). \label{rLpLma}
\end{align}
The half-currents have decompositions
\begin{align*}
   x^{(a)+}(u)&=\sum_{n\ge0}x^{(a)}_n u^{-n-1}, &   x^{(a)-}(u)&=-\sum_{n<0}x^{(a)}_n u^{-n-1},
\end{align*}
where $x^{(a)}_n=(x\otimes z^n,0,0)$ for $n\in\mathbb Z$, $x\in\{h,e,f\}$. This means that this
algebra coincides with a classical limit of the central extension of the Yangian double
$DY(\widehat{\slt})$~\cite{Kh}.

\bb{Case (b): $\tau\to i\infty$ ($\omega=1$, $\omega'\to\infty$, $\tau=\omega'/\omega$, $\Im\tau>0$).}
In this case the degenerations of elliptic Green distributions
look as follows:
\begin{align*}
 G(u,z)&\to\psi^+(u-z), \\
 \G\pm(u,z)&\to\pi\ctg\pi\lambda+\psi^\pm(u,z),
\end{align*}
where
\begin{align*}
\la\psi^\pm(u,z),s(u)\ra_u=\oint\limits_{\substack{|u|>|z| \\ |u|<|z|}}\frac{du}{2\pi i}\pi\ctg\pi(u-z)s(u).
\end{align*}
By the same reason the degenerated algebras are isomorphic for
different $\lambda$ and we shall consider this bialgebra only in
the limit $\lambda\to-i\infty$. The $r$-matrix, $L$-operators and
$rLL$-relations in this case take the form
\begin{align}
r^{(b)+}(u,v)
 =\begin{pmatrix}
   \frac12\psi^+(u,v) & 0              & 0              & 0              \\
  0              & -\frac12\psi^+(u,v) & -\pi i+\psi^+(u,v)      & 0              \\
  0              & \pi i+\psi^+(u,v)       & -\frac12\psi^+(u,v) & 0              \\
  0              & 0              & 0              &  \frac12\psi^+(u,v)
 \end{pmatrix}. \label{rb}
\end{align}
\begin{align*}
 \L^{(b)\pm}(u)&=\lim\limits_{\lambda\to-i\infty}\lim\limits_{\tau\to i\infty}\L_{\lambda}^\pm(u)=
 \begin{pmatrix}
  \frac12h^{(b)\pm}(u) & f^{(b)\pm}(u) \\
  e^{(b)\pm}(u)       & -\frac12h^{(b)\pm}(u)
 \end{pmatrix},
\end{align*}
\begin{align*}
 [\L_1^{(b),\pm}(u),\L_2^{(b)\pm}(v)]&=[\L_1^{(b)\pm}(u)+\L_2^{(b)\pm}(v),r^{(b)+}(u,v)], \\
 [\L_1^{(b)+}(u),\L_2^{(b)-}(v)]&=[\L_1^{(b)+}(u)+\L_2^{(b)-}(v),r^{(b)+}(u,v)]
 +c\cdot\frac{\partial}{\partial u}r^{(b)+}(u,v)
\end{align*}
The half-currents have the following decompositions
\begin{align*}
  h^{(b)+}(u)&=\sum_{n\ge0}h^{(b)}_n \dfrac{\partial^n}{\partial u^n}\ctg\pi u, &
  h^{(b)-}(u)&=-\sum_{n\ge0}h^{(b)}_{-n-1} u^n, \\
  e^{(b)+}(u)&=ie^{(b)}_0+\sum_{n\ge0}e^{(b)}_n \dfrac{\partial^n}{\partial u^n}\ctg\pi u, &
  e^{(b)-}(u)&=ie^{(b)}_0-\sum_{n\ge0}e^{(b)}_{-n-1} u^n, \\
  f^{(b)+}(u)&=-if^{(b)}_0+\sum_{n\ge0}f^{(b)}_n \dfrac{\partial^n}{\partial u^n}\ctg\pi u, &
  f^{(b)-}(u)&=-if^{(b)}_0-\sum_{n\ge0}f^{(b)}_{-n-1} u^n,
\end{align*}
where $x^{(b)}_n=(x\otimes\pi\dfrac{(-1)^n}{n!}z^n,0,0)$,
$x^{(b)}_{-n-1}=(x\otimes\pi\dfrac{(-1)^n}{n!}\dfrac{\partial^n}{\partial z^n}\ctg\pi z,0,0)$
for $n\ge0$, $x\in\{h,e,f\}$.

\subsection{Degeneration of the quasi-Lie bialgebra $\utau$}
\label{sec42}

In the case of algebra $\utau$ there are three cases of degenerations: {\bf(a)}, {\bf(b)} and {\bf(c)}.
The rational degeneration {\bf(a)} and trigonometric degeneration {\bf(b)} are analogous to the
corresponding degenerations of $\cEF$. Additionally there is one more trigonometric case {\bf(c)},
when $\omega\to\infty$ and $\omega'$ is finite. It is not equivalent to the case {\bf(b)} because
the integration contour for $\utau$ is not symmetric in this case.
In the cases {\bf(a)} and {\bf(c)}
the degeneration of elliptic Green distributions  acts on another test function algebra $Z$.
This is an algebra of entire functions $s(u)$ subjected to the inequalities
$|u^n s(u)|<C_n e^{p|\Im u|}$, $n\in\mathbb Z_+$, for some constants $C_n,p>0$
depending on $s(u)$~\cite{GSh1}. The scalar product in $Z$ is $\la s(u),
t(u)\ra_u=\int_{-\infty}^{+\infty}\frac{du}{2\pi i}s(u)t(u)$. The distributions acting on
$K$ can be considered as periodic distributions acting on $Z$.

\bb{Case (a): $\omega\to\infty$, $\omega'\to\infty$, ($\Im\frac{\omega'}{\omega}>0$).}
The degenerating of the elliptic Green distributions in this case reads as
\begin{align*}
 \frac1{\omega}{\cal G}\big(\frac{u-z}{\omega}\big)&\to\frac1{u-z-i0}=\Ph+(u-z), \\
 \frac1{\omega}{\cal G}_{\frac\lambda{\omega}}^\pm\big(\frac{u-z}{\omega}\big)&\to
   \frac{u-z+\lambda}{(u-z\mp i0)\lambda}=\frac1{\lambda}+\Ph\pm(u-z),
\end{align*}
where we introduced the rational Green distributions acting on the test function algebra $Z$ by the formula
\begin{align*}
\la\Ph\pm(u-z),s(u)\ra=\int\limits_{\substack{\Im u<\Im z \\ \Im u>\Im z}}\frac{du}{2\pi i}\frac1{u-z}s(u),
\end{align*}
with infinite horizontal integration lines. They can be represented as integrals
\begin{align*}
 \Ph+(u-z)&=2\pi i\int\limits_{-\infty}^0e^{2\pi i k(u-z)}dk=2\pi i\int
 \limits_0^{+\infty}e^{-2\pi i k(u-z)}dk, \\
 \Ph-(u-z)&=-2\pi i\int\limits_0^{+\infty}e^{2\pi i k(u-z)}dk=-2\pi i\int
 \limits_{-\infty}^0e^{-2\pi i k(u-z)}dk, 
\end{align*}
These formulae are degenerations of both expansions~\eqref{Ggpmdec} and \eqref{GexpanJ}.

As above all the algebras that are obtained from the limit $\omega,\omega'\to\infty$
are isomorphic for the different values of the parameter $\lambda$ and it is sufficient
to describe the limit case $\lambda=\infty$. The $rLL$-relations are the
same as~\eqref{rLpmLpma},~\eqref{rLpLma}
\begin{align*}
 [\LL_{1}^{(a),\pm}(u),\LL_{2}^{(a)\pm}(v)]&=[\LL_{1}^{(a)\pm}(u)+\LL_{2}^{(a)\pm}(v),\rr^{(a)+}(u-v)], \\
 [\LL_{1}^{(a)+}(u),\LL_{2}^{(a)-}(v)]&=[\LL_{1}^{(a)+}(u)+\LL_{2}^{(a)-}(v),\rr^{(a)+}(u-v)]
 +c\cdot\frac{\partial}{\partial u}\rr^{(a)+}(u-v)
\end{align*}
with similar $r$-matrix
\begin{align}
\rr^{(a)+}(u-v)&=\lim_{\lambda\to\infty}\lim_{\omega,\omega'\to\infty}\frac1{\omega}
\rr_{\frac\lambda{\omega}}^+\Big(\frac{u-v}{\omega}\Big)= \notag \\
 &=\begin{pmatrix}
   \frac12\Ph+(u-v) & 0              & 0              & 0              \\
  0              & -\frac12\Ph+(u-v) & \Ph+(u-v)      & 0              \\
  0              & \Ph+(u-v)       & -\frac12\Ph+(u-v) & 0              \\
  0              & 0              & 0              & \frac12\Ph+(u-v)
 \end{pmatrix}, \label{raJinf}
\end{align}
but the entries of the $L$-matrix
\begin{align*}
 \LL^{(a)\pm}(u)&=\lim_{\lambda\to\infty}\lim\limits_{\omega,\omega'\to\infty}
 \dfrac1{\omega}\LL_{\frac{\lambda}{\omega}}^\pm\Big(\dfrac{u}{\omega}\Big)=
 \begin{pmatrix}
  \frac12h^{(a)\pm}(u) & f^{(a)\pm}(u) \\
  e^{(a)\pm}(u)       & -\frac12h^{(a)\pm}(u)
 \end{pmatrix}
\end{align*}
are decomposed to the integrals instead of the series:
\begin{align}
  x^{(a)+}(u)&=\int\limits_0^{+\infty}x^{(a)}_k e^{-2\pi i k u}dk, &
  x^{(a)-}(u)&=-\int\limits_{-\infty}^0x^{(a)}_k e^{-2\pi i k u}dk,
\end{align}
where $x^{(a)}_k=(x\otimes 2\pi i e^{2\pi i k z},0,0)$, $x\in\{h,e,f\}$.
These half-currents form a quasi-classical degeneration of the algebra
${\cal A}_\hbar(\widehat{\slt})$~\cite{KLP-AMS}. The difference between
algebras $DY(\widehat{\slt})$ and ${\cal A}_\hbar(\widehat{\slt})$
is considered in details on the quantum level in this paper.

\bb{Case (b). $\tau\to i\infty$, ($\omega=1$, $\omega'\to\infty$, $\tau=\omega'/\omega$, $\Im\tau>0$).}
Taking the limit $\tau\to i\infty$ in the formulae~\eqref{Ggpmdec} and \eqref{GexpanJ} we obtain
\begin{align*}
 {\cal G}(u-z)&\to\pi i+2\pi i\sum_{n>0}e^{-2\pi in(u-z)}=\tilde\psi^+(u-z), \\
 \Gg\pm(u-z)&\to\pi\ctg\pi\lambda-\pi i
    \pm2\pi i\sum_{\substack{n\ge0 \\ n<0}}e^{-2\pi in(u-z)}=\pi\ctg\pi\lambda+\tilde\psi^\pm(u-z).
\end{align*}
where
\begin{align*}
\la\tilde\psi^\pm(u-z),s(u)\ra=\int\limits_{\substack{\Im u<\Im z \\ \Im u>\Im z}}
\frac{du}{2\pi i}\pi\ctg\pi(u-z)s(u), 
\end{align*}
where $s\in K$ and the integration is taken over a horizontal line segments with unit length.
In these notations the $r$-matrix (in the limit $\lambda\to-i\infty$) can be written as
\begin{align}
&\rr^{(b)+}(u-v)=\lim_{\lambda\to-i\infty}\lim_{\tau\to i\infty}\rr_\lambda^+(u-v)= \notag \\ 
 &=\begin{pmatrix}
   \frac12\tilde\psi^+(u-v) & 0              & 0              & 0              \\
  0              & -\frac12\tilde\psi^+(u-v) & -\pi i+\tilde\psi^+(u-v)      & 0              \\
  0              & \pi i+\tilde\psi^+(u-v)   & -\frac12\tilde\psi^+(u-v) & 0              \\
  0              & 0              & 0              & \frac12\tilde\psi^+(u-v)
 \end{pmatrix}. \label{rbJ}
\end{align}
Setting
\begin{align*}
 \LL^{(b)\pm}(u)&=\lim_{\lambda\to-i\infty}\lim\limits_{\tau\to i\infty}\LL_{\lambda}^\pm(u)=
 \begin{pmatrix}
  \frac12h^{(b)\pm}(u) & f^{(b)\pm}(u) \\
  e^{(b)\pm}(u)       & -\frac12h^{(b)\pm}(u)
 \end{pmatrix}.
\end{align*}
one derives
\begin{align*}
 [\LL_{1}^{(b),\pm}(u),\LL_{2}^{(b)\pm}(v)]=&[\LL_{1}^{(b)\pm}(u)+\LL_{2}^{(b)\pm}(v),\rr^{(b)+}(u-v)], \\
 [\LL_{1}^{(b)+}(u),\LL_{2}^{(b)-}(v)]=&[\LL_{1}^{(b)+}(u)+\LL_{2}^{(b)-}(v),\rr^{(b)+}(u-v)]
 +c\cdot\frac{\partial}{\partial u}\rr^{(b)+}(u-v),
\end{align*}
which define some Lie algebra together with half-current decompositions
\begin{align*}
 h^{(b)+}(u)&=-\frac12h^{(a)}_0+\sum_{n\ge0}h^{(a)}_n e^{-2\pi inu}, &
 h^{(b)-}(u)&=-\frac12h^{(a)}_0-\sum_{n<0}h^{(a)}_n e^{-2\pi inu}, \\
 e^{(b)+}(u)&=\sum_{n\ge0}e^{(a)}_n e^{-2\pi inu}, &
 e^{(b)-}(u)&=-\sum_{n<0}e^{(a)}_n e^{-2\pi inu},\\
 f^{(b)+}(u)&=\sum_{n>0}f^{(a)}_n e^{-2\pi inu}, &
 f^{(b)-}(u)&=-\sum_{n\le0}f^{(a)}_n e^{-2\pi inu}.
\end{align*}
where $x^{(b)}_k=(x\otimes 2\pi i e^{2\pi inz},0,0)$, $x\in\{h,e,f\}$. This is exactly
an affine Lie algebra $\widehat{\slt}$ with a bialgebra structure inherited from the quantum
affine algebra $U_q(\widehat{\slt})$.

\bb{Case (c): $\omega=\to\infty$, $\omega'=\tau\omega=\frac{i}\eta=const$, ($\Re\eta>0$).}
Substituting $u\to\frac{u}{\omega}$, $z\to\frac{z}{\omega}$ to the expansion for ${\cal G}(u-z)$
we yield the following degeneration
\begin{align}
 \frac1\omega{\cal G}\big(\frac{u-z}\omega\big)&\to
  \Psi(u-z)\stackrel{df}{=}2\pi i\vpint\limits_{-\infty}^{+\infty}
   \frac{e^{-2\pi i k(u-z)}dk}{1-e^{-\frac{2\pi k}\eta}}. \label{GdaJ}
\end{align}
Here $\vpints$ means an integral in sense of principal value. The integral in the
formula~\eqref{GdaJ} converges in the domain $\Re\eta^{-1}<\Im(u-z)<0$ and is equal
to $\pi\eta\cth\pi\eta(u-z)$ in this domain. It means that the distribution $\Psi(u-z)$
defined by formula~\eqref{GdaJ} acts on $Z$ as follows
\begin{align}
\la\Psi(u-z),s(u)\ra=\int\limits_{-\Re\eta^{-1}<\Im (u-z)<0}\frac{du}{2\pi i}
\pi\eta\cth\pi\eta(u-z)s(u).\label{Psitg}
\end{align}

The degeneration of Green distributions parameterized by $\lambda$ can be performed in different ways.
We can consider a more general substitution $\lambda\to\mu+\frac\lambda\omega$ instead of
$\lambda\to\frac\lambda\omega$ used above. Substituting $u\to\frac{u}{\omega}$,
$z\to\frac{z}{\omega}$, $\lambda\to\mu+\frac{\lambda}{\omega}$ to the formulae~\eqref{Ggpmdec},
\eqref{GexpanJ} and taking the limits $\omega\to\infty$ and $\lambda\to\infty$  we obtain
\begin{align}
 \frac1{\omega}{\cal G}_{\mu+\frac\lambda{\omega}}^+\big(\frac{u-z}{\omega}\big)&\to
 2\pi i\int\limits_{-\infty}^{\infty}\frac{e^{-2\pi i k(u-z)}dk}{1-e^{-\frac{2\pi k}\eta-2\pi i\mu}},
 \label{GpdaJmu}\\
 \frac1{\omega}{\cal G}_{\mu+\frac\lambda{\omega}}^-\big(\frac{u-z}{\omega}\big)&\to
 2\pi i\int\limits_{-\infty}^{\infty}\frac{e^{-2\pi i k(u-z)}dk}{e^{\frac{2\pi k}\eta+2\pi i\mu}-1}.
 \label{GmdaJmu}
\end{align}
Left hand sides of~\eqref{GpdaJmu} and \eqref{GmdaJmu} as well as right hand sides are invariant under
$\mu\to\mu+1$, but the right hand side is not holomorphic with respect to $\mu$ because the integrand has poles. The complex plane split up to the following  analyticity zones 
$\frac{\Im\eta\Im\mu}{\Re\eta}+n<\Re\mu<\frac{\Im\eta\Im\mu}{\Re\eta}+n+1$, $n\in\mathbb Z$,
and due to periodicity with respect to $\mu$ one can consider only one of these zones.

Integrals in the formulae~\eqref{GpdaJmu} and \eqref{GmdaJmu} converge in the domain
$\Re\eta^{-1}<\Im(u-z)<0$ and $0<\Im (u-z)<\Re\eta^{-1}$ respectively and they can be
calculated like the integral in~\eqref{GdaJ} for the chosen zone. Denote by $\Psi_\mu^+(u-z)$
and $\Psi_\mu^-(u-z)$ the analytic continuation with respect to $\mu$ of the right hand sides
of~\eqref{GpdaJmu} and \eqref{GmdaJmu} respectively from the zone
\begin{align}
 \frac{\Im\eta\Im\mu}{\Re\eta}<\Re\mu<\frac{\Im\eta\Im\mu}{\Re\eta}+1. \label{zone_mu}
\end{align}
Thus, this degeneration of Green distributions can be rewritten as
\begin{align}
 \lim_{\omega\to\infty}{\cal G}_{\mu+\frac\lambda{\omega}}^\pm\big(\frac{u-z}{\omega}\big)
    &=\Psi_\mu^\pm(u-z), \notag \\
\la\Psi_\mu^\pm(u-z),s(u)\ra&=\int\limits_{\substack{-\Re\eta^{-1}<\Im (u-z)<0 \\
0<\Im(u-z)<\Re\eta^{-1}}}\frac{du}{2\pi i}2\pi\eta\frac{e^{-2\pi\eta\mu(u-z)}}{1-e^{-2\pi\eta(u-z)}}s(u)\ ,
 \label{Psipsh}
\end{align}
where $s\in Z$ and the integrals are taken over the horizontal lines.

For the values $\Re\mu=\frac{\Im\eta\Im\mu}{\Re\eta}+n$, $n\in\mathbb Z$, integrands
in~\eqref{GpdaJmu}, \eqref{GmdaJmu} have a pole on the real axis and the distributions
$\Psi_\mu^+(u-z)$ and $\Psi_\mu^-(u-z)$ regularize these integrals as analytical continuation
(see~\cite{GSh1}). The $r$-matrix obtained by another regularization does not satisfy the CYBE.

The degeneration of $r$-matrix is~\footnote{Let us remark that the degeneration of
the entry $r^+(u-v)_{12,21}={\cal G}^+_{-\lambda}(u-v)$ for the zone~\eqref{zone_mu} is
$\Psi^+_{1-\mu}(u-v)$, but is not $\Psi^+_{-\mu}(u-v)$ as one could expect, because of
periodicity with respect to $\mu$ and the fact that $\mu$ belongs to the zone~\eqref{zone_mu}
if and only if $1-\mu$ belongs to the zone~\eqref{zone_mu}. One can also use the relations
${\cal G}^+_{-\mu-\frac\lambda\omega}(\frac{u-v}\omega)=-{\cal G}^-_{\mu+\frac\lambda\omega}
(\frac{v-u}\omega)\longrightarrow-\Psi^-_{\mu}(v-u)=\Psi^+_{1-\mu}(u-v)$.}
\begin{align}
&\rr^{(c)+}(u-v)=\lim_{\lambda\to\infty}\lim_{\omega\to\infty}\frac1{\omega}
\rr_{\mu+\frac\lambda{\omega}}^+\Big(\frac{u-v}{\omega}\Big)= \notag \\ 
 &=\begin{pmatrix}
   \frac12\Psi(u-v) & 0              & 0              & 0              \\
  0              & -\frac12\Psi(u-v) & -\Psi_\mu^-(v-u)      & 0              \\
  0              & \Psi_\mu^+(u-v)   & -\frac12\Psi(u-v) & 0              \\
  0              & 0              & 0              & \frac12\Psi(u-v)
 \end{pmatrix}. \label{rcJ}
\end{align}
The $L$-operators
\begin{align*}
 \LL^{(c)\pm}(u)&=\lim_{\lambda\to\infty}\lim\limits_{\omega\to\infty}\dfrac1
 {\omega}\LL_{\mu+\frac{\lambda}{\omega}}^\pm\Big(\dfrac{u}{\omega}\Big)=
 \begin{pmatrix}
  \frac12h^{(c)\pm}(u) & f^{(c)\pm}(u) \\
  e^{(c)\pm}(u)       & -\frac12h^{(c)\pm}(u)
 \end{pmatrix}.
\end{align*}
with this $r$-matrix satisfy the dynamical $rLL$-relations
\begin{align}
 [\LL_{1}^{(c)\pm}(u),\LL_{2}^{(c)\pm}(v)]=[\LL_{1}^{(c)\pm}(u)&+\LL_{2}^{(c)\pm}(v),\rr^{(c)+}(u-v)]
 -c\cdot i\eta^2\frac{\partial}{\partial\eta}\rr^{(c)+}(u-v)\ , \nn \\
 [\LL_{1}^{(c)+}(u),\LL_{2}^{(c)-}(v)]=[\LL_{1}^{(c)+}(u)&+\LL_{2}^{(c)-}(v),\rr^{(c)+}(u-v)]+
  \label{rLLcJ} \\
 &+c\cdot\bigg(\frac{\partial}{\partial u}-i\eta^2\frac{\partial}{\partial\eta}\bigg)\rr^{(c)+}(u-v)\ . \nn
\end{align}
Decompositions of the half-current in this degeneration are
\begin{align*}
 h^{(c)+}(u)&=\vpint\limits_{-\infty}^{+\infty}h^{(c)}_k
 \frac{e^{-2\pi i k u}dk}{1-e^{-\frac{2\pi k}\eta}}\ , &
 h^{(c)-}(u)&=\vpint\limits_{-\infty}^{+\infty}h^{(c)}_k
 \frac{e^{-2\pi i k u}dk}{e^{\frac{2\pi k}\eta}-1}\ , \\
 e^{(c)+}(u)&=\vpint\limits_{-\infty}^{+\infty}e^{(c)}_k
 \frac{e^{-2\pi i k u}dk}{1-e^{-\frac{2\pi k}\eta-2\pi i\mu}}\ , &
 e^{(c)-}(u)&=\vpint\limits_{-\infty}^{+\infty}e^{(c)}_k
 \frac{e^{-2\pi i k u}dk}{e^{\frac{2\pi k}\eta+2\pi i\mu}-1}\ , \\
 f^{(c)+}(u)&=\vpint\limits_{-\infty}^{+\infty}f^{(c)}_k
 \frac{e^{-2\pi i k u}dk}{1-e^{-\frac{2\pi k}\eta+2\pi i\mu}}\ , &
 f^{(c)-}(u)&=\vpint\limits_{-\infty}^{+\infty}f^{(c)}_k
 \frac{e^{-2\pi i k u}dk}{e^{\frac{2\pi k}\eta-2\pi i\mu}-1}\ ,
\end{align*}
where $x^{(a)}_k=x\otimes 2\pi i e^{2\pi i k z}$, $x\in\{h,e,f\}$.
We do not make explicit the dependence of the parameter $\mu$ because, contrary
to $\lambda$, it is not a dynamical parameter. We also omit dependence on the parameter $\eta$
which provides the dynamics over $c$ just as we omitted its analogue $\tau$ in the elliptic case.
The case $\mu=\frac12$, $\Im\eta=0$ coincides with the quasi-classical limit of the
quantum current algebra ${\cal A}_{\hbar,\eta}(\widehat{\slt})$ \cite{KLP98,CKP}.
This classical algebra was investigated in~\cite{KhLPST} in details. Other degenerations
{\bf(c)} are still not investigated, though the matrices $\rr^{(c)+}(u)$ fit the Belavin-Drinfeld
classification~\cite{BD}.

\section{Averaging of $r$-matrices}
\label{sec5}

Now we will use the averaging method of Faddeev-Reshetikhin~\cite{RF} and will write down trigonometric and elliptic $r$-matrices starting with rational and trigonometric solutions of the Classical Yang-Baxter Equation, respectively. Thereby we show that the $r$-matrix satisfying to a Dynamical Classical Yang-Baxter Equation can be also obtained by this method.

\setcounter{bbcount}{0}

\bb{CYBE.} A meromorphic $\mathfrak{a}\otimes\mathfrak{a}$-valued function $X(u)$
(in our case $\mathfrak a=\slt$) is called solution of the CYBE if it satisfies the equation
\begin{multline}
  [X_{12}(u_1-u_2),X_{13}(u_1-u_3)]+[X_{12}(u_1-u_2),X_{23}(u_2-u_3)]+ \\
  +[X_{13}(u_1-u_3),X_{23}(u_2-u_3)]=0. \label{CYBE}
\end{multline}
The $r$-matrices $\rr^{(a)+}(u-v)$, $\rr^{(b)+}(u-v)$, $\rr^{(c)+}(u-v)$
defined by formulae~\eqref{raJinf}, \eqref{rbJ} and \eqref{rcJ} satisfy CYBE~\eqref{CYBE},
what follows from the fact that they are regularization of the corresponding rational
and trigonometric solutions of CYBE in the domain $\Im u<\Im v$. Indeed, in order to
check the equation~\eqref{CYBE} for these $r$-matrices it is sufficient to check~\eqref{CYBE}
in the domain $\Im u_1<\Im u_2<\Im u_3$. The regularization of the first
two solutions of CYBE ({\bf(a)} and {\bf(b)} cases) but in domain $|u|>|v|$ are
$r$-matrices $r^{(a)+}(u,v)$ and $r^{(b)+}(u,v)$ (formulae~\eqref{ra} and \eqref{rb})
respectively. Hence they also satisfy~\eqref{CYBE}, but where $X_{ij}(u_i-u_j)$ replaced
by $X_{ij}(u_i,u_j)$. The elliptic $r$-matrix $\rr^+_\lambda(u-v)$ satisfies
{\itshape Dynamical} CYBE, but it can be also obtained by the averaging method.

\bb{Basis of averaging.} As it was shown in~\cite{BD} each solution of CYBE $X(u)$ is a
rational, trigonometric or elliptic (doubly periodic) function of $u$, the poles of $X(u)$
form a lattice $\mathfrak R\subset\mathbb C$ and there is a group homomorphism
$A\colon\mathfrak R\to\Aut\mathfrak{g}$ such that for each $\gamma\in\mathfrak R$
one has the relation $X(u+\gamma)=(A_\gamma\otimes\id)X(u)$. Having a rational solution
$X(u)$, for which $\mathfrak R=\{0\}$, and choosing an appropriate automorphisms
$A=A_{\gamma_0}$ we can construct the trigonometric solution with $\mathfrak R=\gamma_0\mathbb Z$
in the form
\begin{align}
 \sum_{n\in\mathbb Z}(A^n\otimes\id)X(u-n\gamma_0). \label{aver}
\end{align}
Applying the same procedure for a trigonometric solution with $\mathfrak R=\gamma_1\mathbb Z$,
 where $\gamma_1/\gamma_0\not\in\mathbb R$, we obtain an elliptic (doubly periodic) solution of
 CYBE with $\mathfrak R=\gamma_1\mathbb Z+\gamma_0\mathbb Z$. The convergence of series in the
 formula~\eqref{aver} should be understood in the principal value sense (below we will detail it).

\bb{Quasi-doubly periodic case.} The entries of elliptic $r$-matrix $\rr^+_\lambda(u)$ --
elliptic Green distributions -- are regularizations of {\itshape quasi-}doubly periodic functions.
This is a direct consequence of the fact that this $r$-matrix satisfies {\itshape Dynamical} CYBE and
therefore does not belong to the Belavin-Drinfeld classification~\cite{BD}. Nevertheless, these
functions have the elliptic type of the pole lattice $\mathfrak R=\Gamma=\mathbb Z+\mathbb Z\tau$
and one can expect that the $r$-matrix $\rr^+_\lambda(u)$ can be represent by formula~\eqref{aver}
with $\gamma_0=\tau$ and $X(u)$ replaced by some trigonometric $r$-matrix with $\mathfrak R=\mathbb Z$.
To pass on from the averaging of meromorphic functions to the averaging of distributions we should
choose the proper regularization. Actually the regularization of this trigonometric $r$-matrix in
this formula can depend on $n$ (see~\eqref{rerb}). The $r$-matrices $r^{(a)+}(u,v)$, $r^{(b)+}(u,v)$,
$r^+_\lambda(u,v)$ can be also regarded as a regularization of the same meromorphic $\slt\otimes\slt$-valued
function, but they depend on~$u$ and~$v$ in more general way than on the difference $(u-v)$.
This makes their averaging
 to be more complicated. By this reason we shall not consider these matrices
 in this section.

\bb{Dynamical elliptic $r$-matrix as an averaging of $\rr^{(b)\pm}(u)$.} To represent
the $r$-matrix $\rr^+_\lambda(u)$ as an averaging of trigonometric matrix~\eqref{rbJ}
we need the following formulae
\begin{align}
 \frac{\theta'(u)}{\theta(u)}&=v.p.\sum_{n\in\mathbb Z}\pi\ctg\pi(u-n\tau)\ , \label{avctg} \\
 \frac{\theta(u+\lambda)}{\theta(u)\theta(\lambda)}&=\frac{\theta'(\lambda)}{\theta(\lambda)}+
 v.p.\sum_{n\in\mathbb Z}\big(\pi e^{-2n\pi i\lambda}\ctg\pi(u-n\tau)
   +(1-\delta_{n0})\pi e^{-2n\pi i\lambda}\ctg\pi n\tau\big)\ , \label{avctg_}
\end{align}
where  $|\Im\lambda|<\Im\tau$, $\lambda\notin\mathbb Z$ and the symbol $v.p.$ means convergence
of the series in the principal value sense:
\begin{align*}
 v.p.\sum_{n\in\mathbb Z}x_n=\lim_{N\to\infty}\sum_{n=-N}^N x_n\ .
\end{align*}
The Fourier expansion of the function $\frac{\theta'(\lambda)}{\theta(\lambda)}$ has
the form~\eqref{GexpanJ} (with $(u-z)$ replaced by $\lambda$) in the domain $-\Im\tau<\Im\lambda<0$.
Substituting this expansion to the right hand side of~\eqref{avctg_} one yields
\begin{align}
 \frac{\theta(u+\lambda)}{\theta(u)\theta(\lambda)}=
  v.p.\sum_{n\in\mathbb Z}\pi e^{-2n\pi i\lambda}\big(\ctg\pi(u-n\tau)+i\big)\ , \label{avctg_p}\\
 \frac{\theta(u-\lambda)}{\theta(u)\theta(-\lambda)}=
  v.p.\sum_{n\in\mathbb Z}\pi e^{2n\pi i\lambda}\big(\ctg\pi(u-n\tau)-i\big)\ .  \label{avctg_m}
\end{align}
The formula~\eqref{avctg_m} is obtained from~\eqref{avctg_p} by replacing $u\to-u$, $n\to-n$,
hence both formulae are valid in the domain $-\Im\tau<\Im\lambda<0$.
Let us choose an automorphism $A=A_\tau$ as follows
\begin{align*}
 &A\colon h \mapsto h, & &A\colon e \mapsto e^{2\pi i\lambda} e,  & &A\colon f \mapsto e^{-2\pi i\lambda} f,
\end{align*}
and define $\vartheta_n=+$ for $n\ge0$ and $\vartheta_n=-$ for $n<0$.
Then the formulae~\eqref{avctg}, \eqref{avctg_p}, \eqref{avctg_m} imply
\begin{align}
 \rr^+_\lambda(u)=v.p.\sum_{n\in\mathbb Z}(A^n\otimes\id)\rr^{(b),\vartheta_n}(u-\tau n)\ ,  \label{rerb}
\end{align}
where $-\Im\tau<\Im\lambda<0$ and $\rr^{(b),-}(u)$ defined by formula~\eqref{rbJ} with
$\tilde\psi^+(u)$ replaced by $\tilde\psi^-(u)$. Let us notice that $r$-matrix $\rr^+_\lambda(u)$
and $\rr^{(b)\pm}(u)$ act as distributions on the same space $K$. Thus we do not have any problem with interpretation of the averaging formula in sense of distributions.

\bb{The matrix $\rr^{(c)+}(u)$ as an averaging of $\rr^{(a)\pm}(u)$.} In this case Green distributions $\Psi^+_\mu(u)$, $\Psi^-_\mu(-u)$, $\Psi(u)$ entering into the $r$-matrix $\rr^{(c)+}(u)$ are defined by~\eqref{Psitg}, \eqref{Psipsh} and the parameter $\mu$ is restricted by~\eqref{zone_mu}. One has the formula
\begin{align*}
 2\pi\eta\frac{e^{2\pi\eta\mu u}}{e^{2\pi\eta u}-1}
  =v.p.\sum_{n\in\mathbb Z}\frac{e^{2\pi i\mu n}}{u-i\eta^{-1}n}. 
\end{align*}
Replacing $u\to-u$, $n\to-n$ in both sides one yields
\begin{align*}
 2\pi\eta\frac{e^{-2\pi\eta\mu u}}{1-e^{-2\pi\eta u}}
  =v.p.\sum_{n\in\mathbb Z}\frac{e^{-2\pi i\mu n}}{u-i\eta^{-1}n}. 
\end{align*}
Let us choose the automorphism $A=A_{i\eta^{-1}}$ in the form
\begin{align*}
 &A\colon H\mapsto H, & &A\colon E\mapsto e^{2\pi i\mu}E, & &A\colon F\mapsto e^{-2\pi i\mu}F\ .
\end{align*}
Then these formulae imply the averaging $r$-matrix
\begin{align*}
 \rr^{(c)+}(u)=\sum_{n\in\mathbb Z}(A^n\otimes\id)\rr^{(a),\vartheta_n}(u-i\eta^{-1}n)\ ,
\end{align*}
where $\frac{\Im\eta\Im\mu}{\Re\eta}<\Re\mu<\frac{\Im\eta\Im\mu}{\Re\eta}+1$ and $\rr^{(a),-}(u)$
is defined by formula~\eqref{raJinf} with $\Phi^+(u)$ substituted by $\Phi^-(u)$. The averaged $r$-matrices $\rr^{(c)+}(u)$ and $\rr^{(a)\pm}(u)$ act also on the same space -- on the
algebra $Z$ from the subsection~\ref{sec42}.

\section*{Acknowledgements} This paper is a part of PhD thesis
of A.S. which he is prepared in co-direction of S.P. and V.R. in
the Bogoliubov  Laboratory of Theoretical Physics, JINR, Dubna and in LAREMA,
D\'epartement de Math\'ematics, Universit\'e d'Angers. He is
grateful to the CNRS-Russia exchange program on mathematical physics  and  personally to J.-M. Maillet
for financial and general support of this thesis project. V.R. are
thankful to B.Enriquez for discussions. He had used during the
project a partial financial support by ANR GIMP, Grant for support of scientific schools
NSh-8065.2006.2
and a support of INFN-RFBR "Einstein" grant (Italy-Russia). He acknowledges a warm hospitality
of Erwin Schr\"{o}dinger Institute for Mathematical Physics and the Program
"Poisson Sigma Models, Lie Algebroids, deformations and higher analogues" where this paper was finished.
S.P. was supported in part by RFBR grant 06-02-17383.

Leonid Vaksman, an excellent mathematician, one of Quantum Group `pioneers', patient teacher and a bright person,
had passed away after coward disease when this paper was finished. We dedicate it to
his memory with a sad and sorrow.
\setcounter{section}{0} \setcounter{subsection}{0}

\renewcommand{\thesection}{\Alph{section}}

\section{Test function algebras $\lfK_0$ and $K=K(\Cyl)$}\label{apA}

\setcounter{bbcount}{0}

\bb{Test function algebra $\lfK_0$.} Let $\lfK_0$ be a set of complex-valued meromorphic
functions defined in some vicinity of origin which have an only pole in the origin. If $s_1(u)$
and $s_2(u)$ are two such functions with domains $U_1$ and $U_2$ then their sum $s_1(u)+s_2(u)$
and their product $s_1(u)\times s_2(u)$ are also functions of this type which are defined
in the intersection $U_1\cap U_2$. Moreover if $s(u)$ is a function from $\lfK_0$
which is not identically zero then there exists a neighborhood $U$ of the origin such that
the domain $U\bs0$ does not contain zeros of function $s(u)$ and, therefore,
the function $\dfrac1{s(u)}$ is a function from $\lfK_0$ with the domain $U$. This
means that the set $\lfK_0$ can be endowed with a structure of a function field. We shall consider
$\lfK_0$ as an associative unital algebra over $\mathbb C$ equipped with the invariant scalar product
\begin{align}
 \la s_1(u),s_2(u)\ra=\oint\limits_{C_0}\frac{du}{2\pi i}s_1(u)s_2(u), \label{lfK_sp}
\end{align}
where $C_0$ is a contour encircling zero and belonging in the intersection of domains
of functions $s_1(u)$, $s_2(u)$, such that the scalar product is a residue in zero. We
consider the algebra $\lfK_0$ as an algebra of test functions. A convergence in $\lfK_0$
is defined as follows: a sequence of functions $\{s_n(u)\}$ converges to
zero if there exists a number $N$ such that all the function $z^Ns_n(u)$ are regular
in origin and all the coefficients in their Laurent expansion tend to zero.
One can consider (instead of the algebra $\lfK_0$ defined in this way one) the completion
$\overline{\lfK}_0=\mathbb C[u^{-1}][[u]]$. Linear continuous functionals on $\lfK_0$
 are called distributions and form the space $\lfK'_0$ (which coincide with $\overline{\lfK}'_0$).
 The scalar product~\eqref{lfK_sp} being continues defines a continuous injection
 $\lfK_0\to\lfK'_0$. We use the notation $\la a(u),s(u)\ra$ for the action of a distribution $a(u)$
 on a test function $s(u)$ and also the notation $\la a(u)\ra_u=\la a(u),1\ra$, where $1$ is a
 function which identically equals to the unit.

One can define a `rescaling' of a test function $s(u)$ as a function $s\big(\frac u\alpha\big)$,
where $\alpha\in\mathbb C$, and therefore a `rescaling' of distributions by the formula
$\la a(\frac u\alpha\big),s(u)\ra=\la a(u),s(\alpha u)\ra$. On the contrary, we are unable to define a
`shift' of test functions by a standard rule, because the operator $s(u)\mapsto s(u+z)$
is not a continuous one~\footnote{Consider, for example, the sum
$s_N(u)=\sum_{n=0}^N(\frac{u}{\alpha})^n$. For each $z$ there exist $\alpha$ such that the
sum $s_N(u+z)$ diverges, when $N\to\infty$.}. Nevertheless we use distributions `shifted' in some sense.
Namely, we say that a two-variable distribution $a(u,z)$ (a linear continuous functional
$a\colon\lfK_0\otimes\lfK_0\to\mathbb C$) is `shifted' if it possesses the properties:
(i) for any $s\in\lfK_0$ the functions $s_1(z)=\la a(u,z),s(u)\ra_u$ and $s_2(u)=\la a(u,z),s(z)\ra_z$
belong to $\lfK_0$; (ii) $\frac{\partial}{\partial u}a(u,z)=-\frac{\partial}{\partial z}a(u,z)$.
Here the subscripts $u$ and $z$ mean the corresponding partial action, for instance,
$\la a(u,z),s(u,z)\ra_u$ is a distribution acting on $\lfK_0$ by the formula
\begin{align*}
 \La\la a(u,z),s(u,z)\ra_u,t(z)\Ra=\la a(u,z),s(u,z)t(z)\ra.
\end{align*}
The condition (ii) means the equality $\la a(u,z),s'(u)t(z)\ra=-\la a(u,z),s(u)t'(z)\ra$.
The condition (i) implies that for any $s\in\lfK_0\otimes\lfK_0$ the expression
\begin{align}
 \la a(u,z),s(u,z)\ra_u=\sum_i\la a(u,z),p_i(u)\ra_u q_i(z), \label{auz_suz}
\end{align}
where $s(u,z)=\sum_i p_i(u)q_i(z)$, belongs to $\lfK_0$ (as a function of $z$).

\bb{Semidirect product.} Now we are able to define a \emph{semidirect product of
two `shifted' distributions} $a(u,z)$ and $b(v,z)$
as a linear continuous functional $a(u,z)b(v,z)$ acting on $s\in\lfK_0\otimes\lfK_0\otimes\lfK_0$
by the rule
\begin{align*}
 \la a(u,z)b(v,z),s(u,v,z)\ra=\La a(u,z),\la b(v,z),s(u,v,z)\ra_v\Ra_{u,z}.
\end{align*}

The `shifted' distribution $a(u,z)$ acting on $\lfK_0\otimes\lfK_0$ can be defined by one
of its partial actions on the function of one variable. For instance, if the partial action
of the type $\la a(u,v),s(u)\ra_u$ is defined for any test function $s(u)$ then one can
calculate the left hand side of~\eqref{auz_suz} and, then, obtain the total action of $a(u,z)$
on the test function $s(u,z)$. This means that a `shifted' distribution (more generally, a
distribution satisfying condition (i)) define a continuous operator on $\lfK_0$.

The main example of a `shifted' distribution is a delta-function $\delta(u,z)$ defined
by one of the formulae
\begin{align*}
 &\la\delta(u,z),s(u)\ra_u=s(z), & &\la\delta(u,z),s(z)\ra_z=s(u), &
 &\la\delta(u,z),s(u,z)\ra_{u,z}=\la s(z,z)\ra_z.
\end{align*}
It is symmetric: $\delta(u,z)=\delta(z,u)$ and one can show that any `shifted'
distribution $a(u,v)$ satisfies
\begin{align}
 a(u,v)\delta(u,z)=a(z,v)\delta(u,z). \label{a_delta}
\end{align}
The distribution $\delta(u,z)$ defines an identical operator on $\lfK_0$.

Another example of a `shifted' distribution is given by Green distributions defined
by~\eqref{Glpmdef}, \eqref{Gdef}. The condition~(i) is obvious if one takes
into account the fact that the functions belonging to $\lfK_0$ have the form
$s(u)\hm=\sum\limits_{n=-N}^\infty s^n\epsilon_{n;\lambda}(u) \hm=\sum\limits_{n=-\infty}^{N}
s_n\epsilon^{n;\lambda}(u)$. The condition~(ii) reads from the formulae~\eqref{Glpmact},
\eqref{Gact}, hence they are actually `shifted' distributions. The condition~(ii) is necessary for
the formula~\eqref{dL} used when one calculates the standard co-cycle $B(\cdot,\cdot)$ on the half-currents.
Indeed, using the partial integration formula
$\la a(u,v)\frac{\partial}{\partial v}x(v)\ra=-\la\frac{\partial}{\partial v}a(u,v)x(v)\ra$, the
formulae for the adjoint action of operator $d$ on total currents~\eqref{xxtc} and
definition~\eqref{hpmht}--\eqref{fpmht} we obtain the formula~\eqref{dL}.

\bb{Test function algebra $K$.} We define the algebra $K$ as an algebra of entire functions
on $\mathbb C$ subjected to the periodicity condition $s(u+1)=s(u)$ and
to the condition $|s(u)|\le C e^{p|\Im u|}$, where the constants $C,p>0$ depend on the function $s(u)$.
The periodicity of these functions means that they can be considered as functions on the
cylinder $\Cyl=\mathbb C/{\mathbb Z}$: $K=K(\Cyl)$. We equip the algebra $K$ with an
invariant scalar product
\begin{align}
 \la s(u),t(u)\ra=\int\limits_{-\frac12+\alpha}^{\frac12+\alpha}\frac{du}{2\pi i}s(u)t(u),
   \qquad s,t\in K, \label{spJ}
\end{align}
which does not depend on a choice of the complex number $\alpha\in\mathbb C$. A
convergence in $K$ is given as follows: a sequence $\{s_n\}\subset K$ tends to zero
if there exist such constants $C,p>0$ that $|s_n(u)|\le C e^{p|\Im u|}$ and for all
$u\in\mathbb C$ the sequence $s_n(u)\to0$. In particular, if $s_n\to0$ then the
functions $s_n(u)$ tends uniformly to zero. Therefore, the scalar product~\eqref{spJ}
is continuous with respect to this topology and it defines a continuous embedding
of $K=K(\Cyl)$ to the space of distributions $K'=K'(\Cyl)$.

Each function $s\in K$ can be restricted to the line segment
$[-\frac12+\alpha;\frac12+\alpha]$, be expanded in this line segment to a Fourier
series and, then, this expansion can be uniquely extend to all the
$\mathbb C$ by the analyticity principle. It means that $\{j_n(u)=e^{2\pi inu}\}_{n\in\mathbb Z}$ is
a basis of $K$ and $\{j^n(u)=2\pi i e^{-2\pi inu}\}_{n\in\mathbb Z}$ is its dual
one with respect to the scalar product~\eqref{spJ}.

The functions belonging to the space $K$ can be correctly shifted because for
all $z\in\mathbb C$ the operator ${\cal T}_z\colon s(u)\mapsto s(u+z)$ is continuous
and maps a periodic function to a periodic one. Hence we can define sifted distributions
in the usual way: $\la a(u-z),s(u)\ra\hm=\la a(u),s(u+z)\ra$. Thereby defined shifted distributions
$a(u-z)$ can be considered as two-variable distributions with properties (i) and (ii) as well as
distributions depending of one of variables as of an argument and of another as of a parameter.
For example the distribution $\delta(u)\in K'$ defined by the formula $\la\delta(u),s(u)\ra=s(0)$
can be shifted by variable $z$ and consider as a distribution of variables $u$ and $z$. This
shifted distribution is called delta-function. Their Fourier expansion looks as follows
\begin{align}
 \delta(u-z)=\sum_{n\in\mathbb Z}j^n(u)j_n(z)=2\pi i\sum_{n\in\mathbb Z} e^{-2\pi in(u-z)}. \label{dffeJ}
\end{align}


\frenchspacing
\begin{thebibliography}{9999999}
{\small

\bibitem[B1]{Bax1} Baxter, R.\,J. Partition function for the eith-vertex lattice model,
{\it Ann. Phys.}, \textbf{70}, (1972), 193--228.

\bibitem[FIJKMY]{FIJKMY}
Foda, O., Iohara, K., Jimbo, M., Kedem, R., Miwa, M., Yan, H.
An elliptic quantum algebra for $\widehat{sl}_2$
{\it  Lett. Math. Phys.} {\bf 32} (1994) 259--268.

\bibitem[L]{L95} Lukyanov, S. Free field representation for
massive integrable models.
{\it Commun. Math. Phys.} {\bf 167} (1995) 183--226.

\bibitem[S]{Skl} Sklyanin E.\,K. Some algebraic structure connected with the Yang-Baxter equation,
Funct. Annal. Appl. \textbf{16}, (1982), 263--270; \textbf{17}, (1983), 273--284.

\bibitem[KLPST]{KhLPST}
Khoroshkin, S., Lebedev, D., Pakuliak, S.,  Stolin, A.,  Tolstoy, V.
{  Classical limit of the scaled eliptic algebra.}
 {\it Compositio Mathematica} {\bf 115} (1999), no. 2, 205--230.


\bibitem[B2]{Bax2} Baxter, R.\,J., Eith-vertex model in lattice statistics and one dimensional
anisotropic Heisenberg chain II: Equivalence to a generalized ice-type lattice model.
{\it Ann. Phys.}, \textbf{76}, (1973), 25--47.


\bibitem[Fe]{F94} Felder, G.  Conformal field theory and
integrable systems associated to
elliptic curves, Proc. ICM Z\"urich 1994, 1247-55, Birkh\"auser (1994);
Elliptic quantum groups,  Proc. ICMP Paris
1994, 211-8, International Press (1995).

\bibitem[ER1]{ER97} Enriquez, B., Rubtsov, V. Quantum groups in higher genus
and Drinfeld's new realizations method ($\mathfrak{sl}_{2}$ case),
{\it Ann. Sci. Ec. Norm. Sup.} {\bf 30} (1997), 821-846.


\bibitem[D1]{D88} Drinfeld, V. New realization of Yangians and quantum
affine algebras. {\it Sov. Math. Dokl.} {\bf 36} (1988) 212--216.

\bibitem[EF]{EF98} Enriquez, B., Felder, G. Elliptic quantum groups
$E_{\tau,\eta}(\mathfrak{sl}_2)$ and quasi-Hopf algebras. {\it Commun. Math. Phys.} {\bf 195}, (1998),
651--689.


\bibitem[D2]{D90} Drinfeld, V. Quasi-Hopf algebras. {\it Leningrad Math. J.}
{\bf 1} (1990) 1419--1457.



\bibitem[BBB]{BBB} Babelon, O., Bernard, D., Billey, E.
A Quasi-Hopf algebra interpretation of quantum 3-j and 6-j symbols and difference equations.
{\it Phys. Lett. B} {\bf 375} (1996) 89--97.

\bibitem[GN]{GN} Gervais, J.-L., Neveu, A. Novel triangle relation and absence of tachyons in
Liouville string field theory. {\it Nucl. Phys.} {\bf 238} (1984) 125.

\bibitem[ER2]{ER-MS} Enriquez, B., Rubtsov, V. Some examples of quantum
groups associated with higher genus algebraic curves, "Moscow Seminar
in Mathematical Physics", 33--65, Amer. Math. Soc. Transl. Ser. 2,
{\bf 191}, Amer. Math. Soc., Providence, RI, 1999.

\bibitem[Fr]{Fr97} Fr\o nsdal, C. Quasi-Hopf deformations of quantum groups.
{\it Lett. Math. Phys.} {\bf 40} (1997) 117--134.




\bibitem[ABRR]{ABRR97} Arnaudon, D., Buffenoir, E., Ragoucy, E., Roche, P.
Universal solutions of quantum dynamical Yang-Baxter equations, {\itshape Lett. Math. Phys.},
\textbf{44}, no. 3, (1998),  201--214.

\bibitem[JKOS1]{JKOS97} Jimbo, M., Konno, H., Odake, S., Shiraishi, J.
Quasi-Hopf twistors for elliptic quantum groups, {\itshape Transform. Groups}, \textbf{4}, no. 4,
(1999), 303--327.


\bibitem[KLP1]{KLP98}
Khoroshkin, S., Lebedev, D.,  Pakuliak, S.
   { Elliptic algebra in the scaling limit.}
{\it Commun. Math. Phys.} {\bf 190} (1998), no. 3, 597--627.

\bibitem[K]{K98} Konno, H., An elliptic algebra $U_{q,p}(\widehat{\slt}_2)$ and the
fusion RSOS model. {\it Commun. Math. Phys.},  \textbf{195}, no. 2, (1998),  373--403.

\bibitem[EPR]{EPR} Enriquez, B., Pakuliak, S., Rubtsov, V., Basic representations of quantum
current algebras in higher genus. Preprint Universit\'e d'Angers UMR 6093, {\tt math/0610398}.


\bibitem[ER3]{ER99} Enriquez, B., Rubtsov, V. Quasi-Hopf algebras
associated with $\mathfrak{sl}_{2}$ and complex curves. {\it Israel Journal
of Mathematics} {\bf 112} (1999) 61--108.


\bibitem[JKOS2]{JKOS} Jimbo, M., Konno, H., Odake, S., Shiraishi, J. Elliptic algebra
$U_{q,p}(\widehat{\mathfrak{sl}}_2)$: Drinfeld currents and vertex operators.
{\it Commun. Math. Phys.} {\bf 199} (1999) 605--647.

\bibitem[Fa]{Fay}
Fay, J.\,D., Theta functions on Riemann surfaces, Springer, (1973).



\bibitem[Kh]{Kh} Khoroshkin, S., Central extension for the Yangian double,
In collection SMF, Colloque ``Septi\`emes Rencontres du Contact Franco-Belge en Alg\`ebre'',
June 1995, Reins.

\bibitem[KLP2]{KLP-AMS}
Khoroshkin, S., Lebedev, D.,  Pakuliak, S.
{\it   Yangian algebras and classical Riemann problem.} "Moscow Seminar in Mathematical
Physics", 163--198,
Amer. Math. Soc. Transl. Ser. 2, {\bf 191}, Amer. Math. Soc., Providence, RI, 1999


\bibitem[LeKP]{CKP}
LeClair, A., Khoroshkin, S., Pakuliak, S.
{Angular quantization of the Sine-Gordon model at the free fermion point.}
{\it Adv. Theor. Math. Phys.} {\bf 3} (1999), no. 5, 1227-1287

\bibitem[RF]{RF}
Reshetikhin, N.\,Yu., Faddeev, L.\,D., Hamiltonian structures for integrable
field theory models, (in Russian)  {\itshape Teoret. Mat. Fiz.}, \textbf{56}, (1983),  no. 3, 323--343.

\bibitem[BD]{BD}
Belavin, A.\,A., Drinfeld, V.\,G., Solutions of the classical Yang-Baxter
 equation for simple Lie algebras, {\itshape Funktsional.
  Anal. i Prilozhen.}, \textbf{16}, (1982), no. 3, 1--29.

\bibitem[WW]{WW}
Whittaker, E.\,T., Watson, G.\,N., A Course of Modern Analysis, 4th eds.,
Cambridge University Press, Cambridge, p. 489 (example 11 in the end of the chapter XXI).

\bibitem[F]{F}
Felder, G., Elliptic quantum groups, Proc. ICMP Paris, 1994, 211-8, International Press, (1995).

\bibitem[GS]{GSh1}
Gelfand, I.\,M., Shilov E.\,G., Generalized functions, Vol.I, Academic Press, New York, (1964).

\bibitem[Vl]{Vl} Vladimirov, V.\,S., Methods Of The Theory Of Generalized Functions,
Taylor and Francis Ltd., 328 pp., (2002).
}\end{thebibliography}
\end{document}